\documentclass{amsart}
\usepackage{amsmath}
\usepackage{amssymb}
\usepackage{amsthm}

\newtheorem{theorem}{Theorem}[section] 
\newtheorem{maintheorem}[theorem]{Main Theorem}
\newtheorem{lemma}[theorem]{Lemma}

\newtheorem{claim}{Claim}[theorem]

\theoremstyle{definition}
 
\newtheorem{definition}[theorem]{Definition}

\theoremstyle{remark}

\newtheorem{remark}[theorem]{Remark}

\newtheorem{discussion}[theorem]{Discussion}

\newcommand{\rest}{\mathord{\restriction}}

\newcommand{\HOM}{{\rm Hom}}
\newcommand{\Ext}{{\rm Ext\, }}
\newcommand{\bP}{{\bf P}}
\newcommand{\bk}{{\bf k}}

\newcommand{\bbz}{{\mathbb Z}}
\newcommand{\acc}{{\rm acc\, }}
\newcommand{\ord}{{\rm ord\, }}
\newcommand{\cf}{{\rm cf}}
\newcommand{\otp}{{\rm otp}}
\newcommand{\vare}{\varepsilon}

\newcommand{\Dom}{{\rm Dom\,}}
\newcommand{\lk}{\langle}
\newcommand{\rk}{\rangle}
\newcommand{\conc}{{}^\frown\!}
\newcommand{\GN}{{\mathfrak N}}

\newcommand{\GB}{{\mathfrak B}}
\newcommand{\GJ}{{\mathfrak J}}
\newcommand{\Gz}{{\mathfrak z}}
\newcommand{\gto}{\stackrel{\rm g}{\to}}
\newcommand{\hto}{\stackrel{\rm h}{\to}}

\begin{document}
\author{Alan Mekler}
\address{Department of Mathematics\\
Simon Fraser University, B.C. Canada}
\thanks{Partially supported by NSERC}

\author{Andrzej Ros{\l}anowski}
\address{Institute of Mathematics\\
 The Hebrew University of Jerusalem\\
 91904 Jerusalem, Israel\\
 and Mathematical Institute of Wroclaw University\\
 50384 Wroclaw, Poland}
\email{roslanow@math.idbsu.edu}
\urladdr{http://math.idbsu.edu/$\sim$roslanow}

\author{Saharon Shelah}
\address{Institute of Mathematics\\
 The Hebrew University of Jerusalem\\
 91904 Jerusalem, Israel\\
 and  Department of Mathematics\\
 Rutgers University\\
 New Brunswick, NJ 08854, USA}
\email{shelah@math.huji.ac.il}
\urladdr{http://www.math.rutgers.edu/$\sim$shelah}
\thanks{Partially supported by the United States --- Israel
Binational Science Foundation. Publication 314.}

\title{On the $p$-rank of $\Ext$}
\date{\today}

\begin{abstract}
\noindent Assume $V=L$ and $\lambda$ is regular smaller than the first
weakly compact cardinal. Under those circumstances and with arbitrary
requirements on the structure of $\Ext(G,\bbz)$ (under well known
limitations), we construct an abelian group $G$ of cardinality $\lambda$
such that for no $G'\subseteq G$, $\vert G'\vert<\lambda$ is $G/G'$ free and
$\Ext(G, \bbz)$ realizes our requirements.
\end{abstract}

\maketitle

\section{Introduction}
In section 2 we give a principle true if ${\bf V}={\bf L}$ which is stronger
than $\diamondsuit^+$ (which was enough for building Kurepa tree); of course
the proof follows proofs of Jensen for diamonds. It seems we do not use its
full strength - we seem to actually need only the cases $\Gz(M^1_\delta)=
\cf(\delta)$ (see \ref{principle}). The principle should be helpful for
building models on $\lambda$ with $\Sigma^1_2$-properties (on $\lambda$).
Also there should be cases where we can prove impossibility (by playing with
those cofinalities). 
\medskip

\noindent In the third section we apply the principle to construct abelian
groups (we drop ``abelian''). For a torsion free group $G$ the group
$\Ext(G,\bbz)$ is divisible and therefore its structure is determined by
ranks $\nu_0(G)$, $\nu_p(G)$ (the numbers of copies of ${\mathbb Q}$ and
$Z(p^\infty)$ in the decomposition of the divisible group, where $p$ ranges
over primes) and $\nu_p(G)\leq 2^{\|G\|}$. By \cite{HHSh91}, if ${\bf
V}={\bf L}$ and a group $G$ is not of the form $G_1\oplus G_2$,
$\|G_1\|<\|G\|$ and $G_2$ free then $\nu_0(G)=2^{\|G\|}$. If $\lambda$ is a
regular cardinal smaller than the first weakly compact cardinal,
$\lambda_p\leq\lambda^+$ (for $p$ prime) then assuming ${\bf V}={\bf L}$ we
construct an abelian torsion free group $G$ such that $\nu_p(G)=\lambda_p$
for each $p$, $\nu_0(G)=\lambda^+$ and $\|G\|=\lambda$.  This result can be
considered as a generalization of a result of Sageev and Shelah which states
the same for $\lambda=\aleph_1$ but under the assumption of CH only (see
\cite{SgSh 138}, for an alternative proof see Eklof and Huber \cite{EH} or
Theorem~XII.2.10 of \cite{EM}).
\medskip

\noindent No advanced knowledge of Group Theory is required, constructions of
the third section are purely combinatorial applications of the principle of
the second section. On the other hand no advance tools of Set Theory are used
--- if one accepts \ref{principle}, the rest is elementary.
\medskip

\noindent {\bf Set theoretical notation}: \hspace{0.15in} For a cardinal
$\kappa$, ${\mathcal H}(\kappa)$ is the family of sets with transitive
closure of cardinality $<\kappa$. If $e$ is a set of ordinals then
$\acc(e)$ denotes the set of accumulation points of $e$ (i.e. limits of
$e$) and $\otp(e)$ is the order type of $e$.
\medskip

\noindent {\bf Group theoretical notation}: \hspace{0.15in} $\bP$ is the set
of prime numbers. As all the groups we shall deal with are abelian, we omit
this adjective. $G, H, K$ denote (abelian) groups, $\oplus$ denotes a direct
sum.  $\bbz$ is the additive group of integers. $\HOM (G, H)$ is the group
of homomorphisms from $G$ to $H$ (with the pointwise addition, i.e.\
$(f+g)(x)=f(x)+g(x)$). If $f\in\HOM(G, \bbz)$ and $p\in\bP$ then $f/p\bbz$
is the following member of $\HOM(G,\bbz/p\bbz)$:
\[(f/p\bbz)(x)=f(x)/p\bbz\quad\mbox{(also called $f(x)+p\bbz$).}\]
For a group $G$ and its subset $\{x_n: n\in I\}\subseteq G$, $\langle x_n:
n\in I\rangle_G$ denotes the subgroup of $G$ generated by $\{x_n: n\in I\}$.
\medskip

\noindent {\bf History}:\hspace{0.15in}  The study of the structure of
$\Ext$ has a long history already. For a review of the main results in the
area we refer the reader to \cite{EM}. 

\noindent The results of this paper were proved in 1986/87 and a preliminary
version of the paper was ready in 1992. The disease and death of the first
author\footnote{Professor Alan H. Mekler died in 1992} stopped his work on
the paper. Later the second author joined in finishing the paper.
\medskip

\noindent {\bf Acknowledgment}:\hspace{0.15in} We would like to thank
Paul~C.~Eklof for helpful comments on earlier versions of the paper. 
\bigskip

\section{A construction principle in ${\bf L}$}

\begin{theorem}[{\bf V}={\bf L}]
\label{principle}
Assume
\begin{enumerate}
\item[(A)] $\lambda$ a regular uncountable not weakly
compact\footnote{can be weaken to the existence of $Y$ (see the
proof)} cardinal.
\item[(B)] $S\subseteq\lambda$ is stationary.
\item[(C)] $\bar P, \bar Q, \bar R$ are finite pairwise disjoint sequences of
predicates and function symbols (so a $\bar P$-model $M$ is $(A,\bar P^M)=
(A,\dots,P_i^M,\dots))$.
\item[(D)] Let $\varphi=\varphi(\bar P,\bar Q,\bar R)$ be a first order
sentence.
\item[(E)] $M^0$ is a $\bar P$-model with universe $\lambda$.
\item[(F)] $E$ is a club of $\lambda$ such that $\delta\in E\Rightarrow
M^0_\delta \stackrel{\rm def}{=} M^0\rest \delta\prec M^0$.
\item[(G)] For $\delta\in E\cup\{\lambda\}$ let $\bk_\delta=\{M^1_\delta:
M^1_\delta$ a $(\bar P\conc \bar Q)$-model expanding $M^0_\delta\}$
and for $M^1_\delta\in\bk_\delta$ let $\bk^+_\delta(M^1_\delta)=\{M^2_\delta:
M^2_\delta$ is a $(\bar P\conc\bar Q\conc\bar R)$-model expanding $M^1_\delta$
and satisfying $\varphi=\varphi(\bar P, \bar Q, \bar R)\}$. Lastly
$\bk^-_\delta=\{M^1_\delta\in\bk_\delta:\bk_\delta^+(M^1_\delta)\neq
\emptyset\}$. 
\end{enumerate}
{\em Then} we can find a well ordering $<^*$ of ${\mathcal H}(\lambda^+)$ of
order type $\lambda^+$, a sequence $\bar e=\lk e_\delta:\delta<\lambda$ is a
singular ordinal $\rk$ and functions $\Gz,\GN,\GB_\vare$ (for
$\vare<\lambda$) such that:
\begin{enumerate}
\item[(a)] The domain of the functions $\Gz$ and $\GN$ is
$\bigcup\limits_{\delta\in E} \bk^-_\delta$; superscript $\delta$ means the
restriction of the function to $\bk_\delta^-$. 
\item[(b)] For $\delta\in E$ and $M^1_\delta\in \bk_\delta^-$ we have:
$\Gz (M^1_\delta)$ is zero or a limit ordinal $<|\delta|^+$.
\item[(c)] For $\delta\in E$, $\GB^\delta_\varepsilon$ is a function with
domain $\{M^1_\delta\in\bk^-_\delta:\Gz(M^1_\delta)>\vare\}$, and
$\GB^\delta=\bigcup\limits_{\vare} \GB^\delta_\vare$.
\item[(d)] For $\delta\in E$, $M^1_\delta\in \bk^-_\delta$ we have
$\GN(M^1_\delta)\in\bk^+_\delta(M^1_\delta)$
\item[(e)] For $\delta \in E$ and $M^1_\delta\in \bk^-_\delta$ we have $\lk
\GB^\delta_\vare(M^1_\delta): \vare<\Gz(M^1_\delta)\rk$ is an increasing
sequence of models, isomorphic to some elementary submodels of $({\mathcal
H}(\lambda^+), \in, <^*)$ but we do not require it to be an elementary chain
nor continuous but we do require the following:
\begin{enumerate}
\item[(*)]\quad $x\in \GB^\delta_{\zeta}(M^1_\delta)\setminus
\GB^\delta_\varepsilon(M^1_\delta), y\in\GB^\delta_\vare(M^1_\delta)
\Rightarrow\GB^\delta_{\zeta}(M^1_\delta)\models y <^*x$.
\end{enumerate}
\item[(f)] For $M^1_\delta\in \bk^-_\delta$, $\vare<\Gz(M^1_\delta)$
the universe of $\GB^\delta_\vare(M^1_\delta)$ is a transitive set
to which $\delta$ belongs.
\item[(g)] For $M^1_\delta\in \bk^-_\delta$  for $\zeta<\Gz(M^1_\delta)$
we have $\lk \GB^\delta_\vare (M^1_\delta):\vare\le \zeta\rk\in
\GB^\delta_{\zeta+1}(M^1_\delta)$ (remember: $\Gz (M^1_\delta)$ if
not zero, is a limit ordinal) and:
\[\GB^\delta_{\zeta+1}(M^1_\delta)\models\mbox{``}
\|\GB^\delta_\zeta(M^1_\delta)\| = \|\delta\|\mbox{''}.\]
\item[(h)] If $M^1_\delta\in \bk^-_\delta$ then  $\GN(M^1_\delta)\notin
\bigcup\{\GB^\delta_\vare(M^1_\delta):\vare<\Gz(M^1_\delta)\}$ and
$M^1_\delta\in\GB^\delta_0(M^1_\delta)$ when $\Gz(M^1_\delta)>0$.
\item[(i)] If $M^1_\lambda\in \bk^-_\lambda$ {\it then} for some
$M^2_\lambda\in \bk^+_\lambda$ we have:
\begin{enumerate}
\item[$(*)$] {\em for some} regular $\sigma\le\lambda$ (possibly $\sigma=0$)
{\em for every} regular $\theta<\lambda$ (the most interesting case is:
$0<\sigma<\lambda\Rightarrow\theta=\sigma$) such that $\{\delta\in
S:\cf(\delta)=\theta\}$ is stationary and {\em for stationarily}\footnote{we
can get a club $E_0\subseteq\lambda$ such that every $\delta\in S\cap E_0$
is OK.} {\em many} $\delta\in S$ we have:
\begin{enumerate}
\item[($\alpha$)] $\cf(\delta)=\theta$
\item[($\beta$)]  $M^2_\lambda\rest\delta\prec M^2_\lambda$
\item[($\gamma$)] $\GN (M^1_\lambda\rest \delta)=M^2_\lambda\rest\delta$.
\item[($\delta$)] $[\sigma<\lambda\Rightarrow \Gz(M^1_\lambda\rest\delta)=
\sigma]$ and $[\sigma=\lambda\Rightarrow \Gz(M^1_\lambda\rest\delta)=
\cf(\delta)]$.
\end{enumerate}
\end{enumerate}
\item[(j)] Suppose $M^1_\lambda\in \bk_\lambda\setminus \bk^-_\lambda$
but $M^1_\lambda\rest\delta\in\bk^-_\delta$ for $\delta\in E$. {\em If}
$M^3_\lambda$ is an expansion of $M^1_\lambda$ with a finite vocabulary
{\em then}\footnote{We can also assign stationary $S'\subseteq S$
were $M^3_\lambda$ was guessed but this is not what we need.}
for some club $E'\subseteq E$ we have
\[\delta\in E'\Rightarrow \{E'\cap \delta,M^3_\lambda\rest\delta\}\in\bigcup
\{\GB^\delta_\vare(M^1_\lambda\rest\delta):\vare<\Gz(M^1_\lambda\rest
\delta)\}.\]
Moreover, {\em if} $M\prec ({\mathcal{H}}(\lambda^+),{\in}, <^*)$ is a proper
$<^*$-initial segment, $\lambda\subseteq M$, $M=\bigcup\limits_{i<\lambda}
M_i$ where $M_i$ increasing continuous, $\|M_i\|<\lambda$ for $i<\lambda$

{\em then} for some club $E'\subseteq E$ for every $\delta\in E'$:
\begin{quotation}
\noindent if ${\bf j}_\delta$ is the Mostowski collapse of $M_\delta$ then

\noindent $\{E'\cap\delta,{\bf j}_\delta[M_\delta]\}\in
\bigcup\limits_{\vare<\Gz(M^1_\lambda\rest\delta)}  
\GB^\delta_\vare(M^1_\lambda\rest\delta)$ and 

\noindent ${\bf j}_\delta[M_\delta]$ is a proper initial segment of
$\bigcup\limits_{\vare<\Gz(M^1_\lambda\rest\delta)} \GB^\delta_\vare
(M^1_\lambda\rest\delta)$ (also the order). 
\end{quotation}
({\em Remark}:\hspace{0.1in} many $M$ has a tower $\langle B^\lambda_\vare:
\vare<\vare^*\rangle$ which collapses to an initial segment.)
\item[(k)] $\bar e$ is a square, i.e. $e_\delta$ a club of $\delta$ of order type $<\delta$ and 
\[\alpha\in \acc(e_\delta)\Rightarrow e_\alpha=e_\delta\cap\alpha,\] 
\[\mbox{if }\ \Gz(M^1_\delta)>0\ \&\ \alpha\in\acc(e_\delta)\ \ \mbox{
then }\ M^1_\alpha\stackrel{\rm def}{=} M^1_\delta\rest\alpha\prec
M^1_\delta\ \mbox{ and }\Gz(M^1_\alpha)>0.\]  
\item[(l)] If $\delta\in E$, $\alpha<\delta$ and $\Gz(M^1_\delta)>0$
then $e_\delta\cap \alpha\in \GB_0^\delta(M^1_\delta)$.
\end{enumerate}
\end{theorem}

\begin{remark}
The interesting case is when the set $S$ satisfies:
\begin{enumerate}
\item[$(\alpha)$] for every $\theta=\cf(\theta)<\lambda$, $\{\delta\in S:
\cf(\delta)=\theta\}$ is stationary,
\item[$(\beta)$] $S$ a set of singular ordinals,
\item[$(\gamma)$] $\lambda=\mu^+\Rightarrow S\subseteq [\mu+1, \lambda)$ and
if $\lambda$ is inaccessible then $S$ is a set of strong limit singular
cardinals and
\item[$(\delta)$] $S$ does not reflect.
\end{enumerate}
\end{remark}

\begin{proof} Let $Y\subseteq \lambda$ be such that for every
$\alpha<\lambda$ the set
\[\{\beta>\alpha: (L_\beta[Y\cap\alpha],\in,Y\cap\alpha)\equiv
(L_{\lambda^+}, \in, Y)\}\]
is bounded in $|\alpha|^+$ (e.g.\ $Y$ is a non reflecting stationary subset of
$\lambda$).

\noindent Let $\bar e=\lk \bar e_\delta:\delta<\lambda$ singular limit
ordinal$\rk$ be as defined by Jensen \cite{Jn}.

\noindent Let $<^*$ be the canonical well ordering of $L$.

\noindent Supose now that $\delta\in E\cup\{\lambda\}$ and $M^1_\delta\in
\bk_\delta$. If $\delta<\lambda$ and $M^1_\delta\in\bk_\delta^-$ then we let
$\GN(M^1_\delta)$ be the $<^*$-first member of $\bk^+_\delta(M^1_\delta)$. Let
\[\begin{array}{ll}
W^1_\delta(M^1_\delta)\stackrel{\rm def}{=}\{\alpha>\delta:&
L_\alpha[M^1_\delta,Y\cap\delta]\cap 
\bk^+_\delta(M^1_\delta)=\emptyset\mbox{ and}\\
\ &(L_\alpha[M^1_\delta,Y\cap \delta], \in, Y\cap\delta)\mbox{ is
elementarily equivalent }\\
\ &\mbox{to }(L_{\lambda^+},\in, Y),\mbox{ moreover it is isomorphic to
some}\\ 
\ &\mbox{elementary submodel of it }(\mbox{demand
}\delta=\sup(Y\cap\delta)\}. 
\end{array}\]
Let $W^2_\delta(M^1_\delta)=\acc(W^1_\delta(M^1_\delta))$. If
$\delta<\lambda$, $W^2_\delta(M^1_\delta)=\emptyset$ then we let
$\Gz(M^1_\delta)=0$. Otherwise 
\[\Gz(M^1_\delta)=\cf\Big(W^1_\delta(M^1_\delta)\cap \sup W^2_\delta
(M^1_\delta)\Big),\]
(so we loose at most finitely many members of $W^1_\delta(M^1_\delta)$) and
let 
\[a_\delta=a_\delta[M^1_\delta]\subseteq W^1_\delta(M^1_\delta)\cap\sup
W^2_\delta(M^1_\delta)\]
be unbounded of order type $\Gz(M^1_\delta)$ and $\gamma\in a_\delta
\Rightarrow a_\delta\cap \gamma\in L_\gamma [M^1_\delta,Y\cap \delta]$ (use
the definition of the square in Jensen \cite{Jn} for $L[Y\cap\delta]$ and
the ordinal $\sup(W^1_\delta(M^1_\delta)\cap\sup W^2_\delta(M^1_\delta)$). 
Let $\GB^\delta_\vare (M^1_\delta)$ be $(L_\alpha[M^1_\delta, Y\cap \delta],
\in,<^*)$ where $\alpha$ is the $\vare$-th member of $a_\delta$.

Let us show that it works, i.e. that the clauses {\bf (a)}--{\bf (l)}
are satisfied.
\begin{description}
\item[Clause (a)] Directly from the choice.
\item[Clause (b)] By the use of $\sup W^2_\delta(M^1_\delta)$.
\item[Clause (c)] Directly from the choice of $\GB^\delta_\vare$ and
$\Gz(M^1_\delta)$.
\item[Clause (d)] By the choice of $\GN(M^1_\delta)$.
\item[Clause (e)] Since $\GB^\delta_\vare(M^1_\delta)$ is
$L_\gamma[M^1_\delta,Y\cap\delta]$ for $\gamma$ the $\vare$-th member of
$a_\delta=a_\delta[M^1_\delta]$ and $a_\delta[M^1_\delta]$ is a subset of
$W^1_\delta(M^1_\delta)$ (see its definition) we are sure that
$\GB^\delta_\vare(M^1_\delta)$ is OK: they are increasing with $\vare$
as the $\vare$-th member of $a_\delta[M^1_\delta]$ increases with
$\vare$; {\bf (*)} is satisfied as $<^*$ is the canonical well ordering
of $L$ (or $L[Y]$, in our model not a big difference).
\item[Clause (f)] See the choice of $W^1_\delta(M^1_\delta)$.
\item[Clause (g)] It follows from $L_\alpha[M^1_\delta, Y\cap\delta]\in
L_{\alpha+1}[M^1_\delta,Y\cap\delta]$, the definition of $\GB^\delta_\vare
(M^1_\delta)$ (and the presence of $a_\delta[M^1_\delta]$), etc.
\item[Clause (h)] It is a consequence of the first clause in the
definition of $W^1_\delta(M^1_\delta)$: $L_\alpha[M^1_\delta,
Y\cap\delta]\cap\bk^+_\delta(M^1_\delta)=\emptyset$.
\item[Clause (i)] Let $M^2_\lambda\in\bk^+_\lambda(M^1_\lambda)$ be
$<^*$-minimal. Define $W^1_\lambda(M^1_\lambda)$ as above. It is bounded
in $\lambda^+$ as $M^2_\lambda\in L_{\lambda^+}$. Then define
$W^2_\lambda(M^1_\lambda)$, $a_\lambda[M^1_\lambda]$ as above and let
$\gamma^*=\sup(W^2_\lambda(M^1_\lambda))$ (so it is a limit ordinal).
Let $\sigma=\cf(\gamma^*)$. Let $\gamma<\lambda^+$ be such that
\[(L_\gamma,{\in},M^1_\lambda,Y,E,S,M^2_\lambda)\prec(L_{\lambda^+},{\in},
M^1_\lambda,Y, E,S,M^2_\lambda).\]
Let $L_\gamma=\bigcup\limits_{i<\lambda}N_i$, $\|N_i\|<\lambda$, $\langle
N_i: i<\lambda\rangle$ increasing continuous and such that
\begin{itemize}
\item $\{E,S,M^1_\lambda, M^2_\lambda,Y,W^1_\lambda(M^1_\lambda),
W^2_\lambda(M^1_\lambda), \gamma^*, a[M^1_\gamma]\}\in N_0$,
\item $N_i\prec (L_\gamma,{\in})$
\item $\langle N_j: j\leq i\rangle \in N_{i+1}$.
\end{itemize}
Let $E'=\{i<\lambda: N_i\cap\lambda=i\}$, $\delta\in S\cap E'$ and let
${\bf j}_\delta$ be the Mostowski collapse of $N_\delta$. Note: ${\bf
j}_\delta(<^*\restriction N_\delta)=<^*_{{\bf j}_\delta[N_\delta]}$, etc.
Now clearly ${\bf j}_\delta$ maps $M^2_\lambda$ to $\GN(M^1_\delta)$,
$W^1_\lambda(M^1_\lambda)$ to $W^1_\delta(M^1_\delta)$,
$W^2_\lambda(M^1_\lambda)$ to $W^2_\delta(M^1_\delta)$ and $a[M^1_\lambda]$
to $a[M^1_\delta]$. If $\sigma<\lambda$ then necessarily
\[\hspace{-0.15in} a[M^1_\delta]={\bf j}_\delta[a[M^1_\lambda]\cap
N_\delta]=\{{\bf j}(\beta): \beta\in a[M^1_\lambda]\cap N_\delta\}=\{{\bf
j}(\beta): \beta\in  a[M^1_\lambda]\}\]
and if $\sigma=\lambda$ then
\[\hspace{-0.15in} a[M^1_\delta]={\bf j}_\delta[a[M^1_\lambda]\cap
N_\delta]= \{{\bf j}_\delta(\beta):\beta\in a[M^1_\lambda], \beta \mbox{ is
} {<}\delta\mbox{-th member of } a[M^1_\lambda]\}.\]
\end{description}
Similarly we can check {\bf (j)}, {\bf (k)}, {\bf (l)}. This finishes the
proof. 
\end{proof}

\begin{remark}
More generally we can phrase parallels of the squared diamond and/or
diamond$^+$.  
\end{remark}

\begin{discussion}
{\em What is the point of this principle?}

\noindent You can just read the next section to see how it works. Still let
us try to explain it. Diamonds on $\lambda$ has been very good in helping to
build a structure $M$ with universe $\lambda$ satisfying some $\Pi^1_2$
statement (like being Souslin).
\smallskip

\noindent We are given $\bar P$. We build here by induction on $\delta\in S$
an increasing sequence of models $M^1_\delta=(\delta,{\bar
P}\restriction\delta, {\bar Q}^\delta)$ carrying some induction
hypothesis. We want to have at the end that there is no ${\bar R}$ such that
$(\lambda,{\bar P}, \bigcup\limits_{\alpha}{\bar Q}^\alpha,{\bar
R})\models\varphi$, so at stage $\delta$ we look at $\GN(M^1_\delta)$ as a
candidate for the bad phenomena. For $\vare<\Gz(M^1_\delta)$ (if
$\Gz(M^1_\delta)=0$ then our life is easier) in
$\GB^\delta_{\vare+1}(M^1_\delta)$ we know
\[\GB^\delta_{\vare+1}(M^1_\delta)\models\mbox{``}\GB^\delta_\vare(M^1_\delta)
\mbox{ is a transitive set of cardinality }\delta\mbox{''}.\]
So we can list all elements of $\GB^\delta_\vare(M^1_\delta)$ in
$\GB^\delta_{\vare+1}(M^1_\delta)$, i.e. we have $f\in\GB^\delta_{\vare+1}
(M^1_\delta)$, $f:\delta\longrightarrow\GB^\delta_\vare(M^1_\delta)$
which is one-to-one and onto. From the outside point of view $\delta$
has small cofinality and by {\bf (k)}+{\bf (l)} one can find a sequence
$\langle\alpha_\xi: \xi<\cf(\delta)\rangle$ cofinal in $\delta$ and such
that every proper initial segment is in $\GB^\delta_{\vare+1}(M^1_\delta)$
(even in $\GB^\delta_0(M^1_\delta)$, usually even in $L_\gamma$). So we
have a fair chance to diagonalize over those sets to fulfill the
obligation in the inductive construction of ${\bar Q}^\delta$, while
``destroying'' the possibility of $\GN(M^1_\delta)$.
\smallskip

\noindent But doing it for one $\vare$ does not suffice. However, if
$\cf(\Gz(M^1_\delta))=\cf(\delta)$ then we can do better. We can find
$f:\delta\longrightarrow\GB^\delta(M^1_\delta)$, one-to-one and onto and
such that $(\forall \alpha<\delta)(f\rest\alpha\in\GB^\delta (M^1_\delta))$
(remember: $\GB^\delta(M^1_\delta)=\bigcup\limits_{\vare<\Gz(M^1_\delta)}
\GB^\delta_\vare(M^1_\delta)$). This is possible as
$\langle\GB^\delta_\vare(M^1_\delta):\vare\leq\zeta\rangle \in
\GB^\delta_{\zeta+1}(M^1_\delta)$, so we can choose
$f_\vare\in\GB^\delta_{\vare+1}(M^1_\delta)$ as the first one-to-one mapping
from $\delta$ onto $\GB^\delta_{\vare+1}(M^1_\delta)$.  So by a demand
$\langle f_\vare: \vare\leq\zeta\rangle\in\GB^\delta_{\vare+1}
(M^1_\delta)$. Now by an easy manipulation we can combine them (using
$\langle\beta_\vare:\vare<\cf(\delta)\rangle$).
\smallskip

\noindent In the proof of \ref{1.2} and \ref{1.2B}, to make
$\nu_p(G)=\lambda_p$, we build together with $G_\alpha$ also $f^{p,\zeta}$
(for $\zeta<\lambda_p$). We need that all non trivial combinations
$\sum\limits_{\ell<n}a_\ell f^{p,\zeta^\ell}\in\HOM(G,\bbz/p\bbz)$ are not
of the form $f/p\bbz$.  This could be a typical application of the
diamond. But we also need that for every $f\in\HOM(G,\bbz/p\bbz)$ there will
be $f'=\sum \limits_{\ell<n}a_\ell f^{p,\zeta^\ell}$ and
$f^*\in\HOM(G,\bbz)$ such that $f-f'=f^*/p\bbz$. For this the normal thing
is to apply $\diamondsuit^+_\lambda$ and to choose $\langle a_\ell:
\ell<n\rangle$ and $f^*\rest G_\delta$ by giving them to approximations of
$f$. But the two demands seem to be hard to go together without what was
said above.
\end{discussion}

\section{Building abelian groups}

One can think of $\Ext(G, K)$ as essentially the family of isomorphism types
of models $(K,H,G,g,h,c)_{c\in K\cup G}$ such that (in our case $K,H,G$ are
abelian groups; we will not mention this usually, and) $h$ is an embedding
from $K$ into $H$, $g$ a homomorphism from $H$ onto $G$ with the range of
$h$ being the kernel of $g$ (i.e. $0\to K{\hto}H{\gto} G\to 0$ being exact)
up to isomorphism over $K\cup G$. Moreover it has a natural additive
structure.  So $\Ext(G, K)=0$ if and only if for any $0\to K{\hto} H{\gto}
G\to 0$ as above, the range of $h$ is a direct summand of $H$.

\noindent We shall not define $\Ext(G,\bbz)$ fully as below we shall quote
theorems characterizing it in a convenient way in the relevant cases. 

\noindent In this section we show how to construct a group $G$ such that
$\Ext(G,\bbz)$ satisfies pre-given requirements (within well-known
limitations, see below). The main tools in the construction are
\ref{principle} and, as a kind of a single step, \ref{1.2B} below.

\begin{definition}
The quantifier $(\forall^* i<\lambda)$ means ``for every large enough
$i<\lambda$'', so this is an abbreviation for ``$(\exists
j<\lambda)(\forall i\in (j,\lambda))$''.
\end{definition}

\begin{definition}
\begin{enumerate}
\item For a sequence $\bar\lambda=\lk \lambda_\ell:\ell<n\rk$ ($n<\omega$)
of pairwise distinct infinite regular cardinals, $I_{\bar \lambda}$ is
the ideal on $\Dom(I_{\bar \lambda})=\prod\limits_{\ell<n}\lambda_\ell$
(called its domain) such that
\[\hspace{-0.3in}A\in I_{\bar \lambda}\quad \mbox{iff}\quad (\forall^*
i_0<\lambda_0)(\forall^* i_1<\lambda_1)\dots(\forall^* i_{n-1}<\lambda_{n-1})
[\lk i_0, \dots, i_{n-1}\rk\notin A].\]
\item For any $\lambda\ge \aleph_0$ we define $\GJ_\lambda$
as the set of ideals of the form $I_{\bar \lambda}$ with $\lambda=
\max\{\lambda_\ell:\ell<n\}$. Let $\GJ_{\le\lambda}=\bigcup\limits_{\mu\le
\lambda} \GJ_{\mu}$.
\end{enumerate}
\end{definition}

\begin{lemma}
\label{decreasing}
Suppose that $\bar{\lambda}=\lk\lambda_\ell: \ell<n\rk$ is a sequence of
pairwise distinct infinite regular cardinals and $\bar{\lambda}^{\rm
dec}=\lk \lambda^\prime_\ell: \ell<n\rk$ is the re-enumeration of
$\bar{\lambda}$ in the decreasing order. Let $\pi:\prod\limits_{\ell<n}
\lambda_\ell\longrightarrow\prod\limits_{\ell<n}\lambda^\prime_\ell$ be the
canonical bijection (i.e.\ $\pi(\eta)(\ell_0)=\eta(\ell_1)$ provided
$\lambda^\prime_{\ell_0}=\lambda_{\ell_1}$). Then
\[A\in I_{\bar{\lambda}}\ \ \ \Rightarrow\ \ \ \pi[A]\in
I_{\bar{\lambda}^{\rm dec}}.\]
\end{lemma}

\begin{proof}
This is an iterated application of the following observation:
\begin{claim}
Let $\lambda_0<\lambda_1$ be regular infinite cardinals, $\psi(x,y,\bar{z})$
be a formula. Then
\[(\forall^* i_0<\lambda_0)(\forall^*
i_1<\lambda_1)\psi(i_0,i_1,\bar{\imath})\ \ \Rightarrow\ \ (\forall^*
i_1<\lambda_1)(\forall^* i_0<\lambda_0)\psi(i_0,i_1,\bar{\imath}).\]
\end{claim}
The claim should be clear and so the lemma.
\end{proof}

\begin{theorem}[V=L]
\label{1.2B}
Assume $\lambda$ is a regular cardinal smaller than the first uncountable
weakly compact cardinal. Suppose that $I_k\in\GJ_{\le\lambda}$, for 
$k<k^* <\omega$, $H$ is a free group with the free basis
\[\{x^k_t:t\in\Dom(I_k)\mbox{ and }k<k^*\}.\]
Further, let $p\in \bP$ and let $f^*\in\HOM(H,\bbz)$ be a homomorphism
such that for some $\ell_0<k^*$
\[\{t\in \Dom(I_{\ell_0}): f^*(x_t^{\ell_0})= 0\}\in I_{\ell_0}.\] 
{\em Then} there is a free group $G$, $H\subseteq G$ such that
$\|G\|=\lambda$, $G/H$ is $\lambda$-free and: 
\begin{enumerate}
\item[$(\alpha)$] there is no $f\in \HOM(G, \bbz)$ extending $f^*$,
\item[$(\beta)$] if $k'< k^*$, $A\in I_{k'}$ then $G/\lk x_t^k: [k=k'\
\&\ t\in A] \mbox{ or } k\neq k'\rk_H$ is free,
\item[$(\gamma)$] if a homomorphism $g^*\in\HOM(H,\bbz)$ is such
that for every $k'<k^*$
\[\{t\in\Dom(I_{k'}): g^*(x^{k'}_t)\neq 0\}\in I_{k'}\]
and $g^+\in\HOM(G,\bbz/p\bbz)$ extends $g^*/p\bbz$ then there exists
$g\in\HOM(G,\bbz)$ such that $g/p\bbz=g^+$ and $g^*\subseteq g$,
\item[$(\delta)$] if $q\in\bP$, $h\in\HOM(H,\bbz/q\bbz)$ is such that for
every $k'<k^*$
\[\{t\in\Dom(I_{k'}): h(x^{k'}_t)\neq 0\}\in I_{k'}\]
then $h$ can be extended to an element of $\HOM(G,\bbz/q\bbz)$.
\end{enumerate}
\end{theorem}

\begin{proof}
Due to lemma \ref{decreasing} it is enough to prove the theorem under the
assumption that the ideals $I_k$ are determined by decreasing sequences
$\bar{\lambda}^k$ of regulars. The proof is by induction on $\lambda$. To
carry out the induction we need the existence of stationary non reflecting
sets and $\diamondsuit^+_\lambda$ only. However we will use this opportunity
to show a simpler application of \ref{principle} and instead of the diamonds
we will use our principle. The construction of \ref{1.2}, though more
complicated, will be similar to the one here.
\medskip

\noindent For $k<k^*$ let $I_k=I_{\bar\lambda^k}$, $\bar \lambda^k=\lk
\lambda^k_\ell:\ell<n_k\rk$, $\lambda_{\ell}^k\leq\lambda$. Thus $\Dom(I_k)
=\prod\limits_{\ell<n_k}\lambda^k_\ell$ and according to what we noted
earlier we assume that the sequences $\bar{\lambda}^k$ are decreasing.
\medskip

\noindent If $\lambda=\aleph_0$ then $\bigwedge\limits_k
n_k=1$:\hspace{0.15in} this case is easy and can be concluded from
\cite{EM}, pp 362--363. However for the sake of the completeness we will sketch
the construction (skipping only some technical details). The following claim
gives us slightly more than needed:
\begin{claim}
Suppose that $H$ is a free group with basis $\{x_n: n\in\omega\}$,
$f^*\in\HOM(H,\bbz)$ is a homomorphism such that $(\forall
N\in\omega)(\exists n>N)(f^*(x_n)\neq 0)$. Then there is a free group
$G\supseteq H$ such that $G/H\cong {\mathbb Q}$ and
\begin{enumerate}
\item there is no $f\in\HOM(G,\bbz)$ extending $f^*$,
\item if $A\subseteq\omega$ is infinite and $h\in\HOM(H,\bbz/q\bbz)$ is such
that 
\[{\rm Ker}(h)\supseteq\langle x_n: n\in A\rangle_G\]
then $G/\langle x_n: n\notin A\rangle_G$ is free and $h$ can be extended to
a homomorphism from $G$ to $\bbz/p\bbz$,
\item if $g^*\in\HOM(H,\bbz)$ is such that $(\exists N\in\omega)(\forall
n>N)(g^*(x_n)=0)$ and $g^+\in \HOM(G,\bbz/p\bbz)$ extends $g^*/p\bbz$\\
then there is $g\in \HOM(G,\bbz)$ such that $g^*\subseteq g$ and
$g/p\bbz=g^+$. 
\end{enumerate}
\end{claim}

\begin{proof}[Proof of the claim]
Let $A_0=\{n\in\omega:f^*(x_n)\neq 0\}$ and let $\{r_n:n\in A_0\}$ enumerate
$\bbz$. Choose inductively positive integers $s_n$ and integers $m_n$ such
that 
\begin{enumerate}
\item[(a)] $m_n\in\{-1,1\}$
\item[(b)] if $n\notin A_0$ then $m_n=1$, $s_n=(n+p)!$
\item[(c)] if $n\in A_0$ then

\noindent $r_n+s_0\cdot\ldots\cdot s_{n-1} m_n f^*(x_n) +
s_0\cdot\ldots\cdot s_{n-2} m_{n-1} f^*(x_{n-1}) + \ldots + s_0m_1 f^*(x_1) +
m_0 f^*(x_0)\neq 0$

\noindent $s_n=(n+p)!\cdot |r_n+s_0\cdot\ldots\cdot s_{n-1} m_n f^*(x_n) +
s_0\cdot\ldots\cdot s_{n-2} m_{n-1} f^*(x_{n-1}) + \ldots + s_0m_1 f^*(x_1) +
m_0 f^*(x_0)|$
\end{enumerate}
Now, let $G$ be the group generated freely by $\{y_n:n\in\omega\}\cup\{x_n:
n\in\omega\}$ except that
\begin{enumerate}
\item[$(*)$] \hspace{0.15in} $s_n y_{n+1}= y_n+m_n x_n$.
\end{enumerate}
Note that the condition $(*)$ implies that for each $k>0$, $n\in\omega$
\begin{enumerate}
\item[$(**)^k_n$] $y_n=s_n s_{n+1}\cdot\ldots\cdot s_{n+k} y_{n+k+1} -
[s_n s_{n+1}\cdot\ldots\cdot s_{n+k-1} m_{n+k} x_{n+k} + 
s_n \cdot\ldots\cdot s_{n+k-2} m_{n+k-1} x_{n+k-1} + \ldots + s_n m_{n+1}
x_{n+1}+ m_n x_n]$.
\end{enumerate}
\medskip

\noindent {\bf 0.}\ \ \ $G$ is freely generated by $\{y_n:n\in\omega\}$ and
$G/H\cong {\mathbb Q}$.
\medskip

\noindent {\bf 1.}\ \ \ There is no $f\in\HOM(G,\bbz)$ extending $f^*$.

\noindent Why? By $(**)^n_0$ the value of $f$ at $y_0$ determines
$f(y_{n+1})$ in the way that is excluded by the choice of $s_n$ for $n\in
A_0$ (clause {\bf (c)}). 
\medskip

\noindent {\bf 2.}\ \ \ If $A\subseteq\omega$ is infinite then $G/\langle
x_n: n\notin A\rangle$ is free.

\noindent Why? Let $\{n_k: k\in\omega\}=A$ be the increasing enumeration.
Let $G_i=\langle y_n: n_{i-1}<n\leq n_i\rangle_G$, $H_i=\langle x_n:
n_{i-1}<n< n_i\rangle_G$ (with a convention that $n_{-1}=-1$). Then
$G=\bigoplus_{i\in\omega} G_i$, $H_i=G_i\cap\langle x_n: n\notin A\rangle_G$
and $\langle x_n: n\notin A\rangle_G=\bigoplus_{i<\omega} H_i$. The groups
$G_i/H_i$ are (freely) generated by $y_{n_i}/H_i$. Hence $G/H$ is free. 
Extending suitable homomorphisms into $\bbz/q\bbz$ should be clear now. 
\medskip

\noindent {\bf 3.}\ \ \ If $g^*\in\HOM(H,\bbz)$ is such that $(\exists
N)(\forall n>N)(g^*(x_n)=0)$ and $g^+\in \HOM(G,\bbz/p\bbz)$ extends
$g^*/p\bbz$ then there is $g\in \HOM(G,\bbz)$ such that $g^*\subseteq g$ and
$g/p\bbz=g^+$. 

\noindent Why? First note that there is at most one homomorphism
$g^+\in\HOM(G,\bbz/p\bbz)$ such that $g^+\supseteq g^*/p\bbz$. This is
because $G/H\cong{\mathbb Q}$:\ \ \ \ if $g^+_1,g^+_2\in\HOM(G,\bbz/p\bbz)$
agree on $H$ then $g^+_1-g^+_2\in\HOM(G,\bbz/p\bbz)$, ${\rm
Ker}(g^+_1-g^+_2)\supseteq H$ and hence
$(g^+_1-g^+_2)/H\in\HOM(G/H,\bbz/p\bbz)$. But the only homomorphism of
${\mathbb Q}$ into $\bbz/p\bbz$ is the trivial one.

Hence it is enough to show that there is an extension $g$ of $g^*$ to a
member of $\HOM(G,\bbz)$ (as then necessarily $g/p\bbz=g^+$ by the
uniqueness). 

Now let $N$ be such that $(\forall n>N)(g^*(x_n)=0)$. Define
\smallskip

\noindent $g(y_n)=0$ for $n>N$,\\
$g(y_N)=-m_N g^*(x_N)$,\\
$g(y_n)=-[s_n s_{n+1}\cdot\ldots\cdot s_{N-1} m_N g^*(x_N)+
s_n\cdot\ldots\cdot s_{N-2}m_{N-1}g^*(x_{N-1})+\ldots+ s_n m_{n+1}
g^*(x_{n+1}) + m_n g^*(x_n)]$ for $n<N$
\smallskip

\noindent and extend it to a homomorphism from $\HOM(G,\bbz)$. Clearly this
$g$ satisfies $g^*\subseteq g$. This finishes the proof of the claim.
\end{proof}

\noindent Assume now that $\lambda$ is smaller than the first weakly compact
uncountable cardinal, $\lambda>\aleph_0$ and below $\lambda$ the theorem
holds.
\medskip

We may think that for some $k_0\leq k^*$ we have
\[k<k_0\ \Rightarrow\ \lambda^k_0<\lambda\ \ \mbox{ and }\ \ k_0\leq k<k^*\
\Rightarrow\ \lambda^k_0=\lambda.\]
Of course we may assume that $k_0<k^*$ (otherwise the inductive hypothesis
applies directly). 

Recall that $\ell_0<k^*$ is such that 
\[\{t\in\Dom(I_{\ell_0}): f^*(x^{\ell_0}_t)= 0\}\in I_{\ell_0}.\] 
Let $\alpha_0$ be defined as follows:
\smallskip

\noindent if $\ell_0<k_0$ then $\alpha_0=0$,

\noindent if $k_0\leq\ell_0<k^*$ and $n_{\ell_0}>1$ then $\alpha_0$ is (the
first) such that 
\[(\forall\alpha>\alpha_0)(\forall^* i_1<\lambda^{\ell_0}_1)\dots (\forall^*
i_{n_{\ell_0}-1}<\lambda^{\ell_0}_{n_{\ell_0}-1})(
f^*(x^{\ell_0}_{\langle\alpha,i_1,\ldots,i_{n_{\ell_0}-1}\rangle})\neq 0),\] 

\noindent if $k_0\leq\ell_0<k^*$ and $n_{\ell_0}=1$ then $\alpha_0<\lambda$ is 
the first ordinal such that $(\forall i_0>\alpha_0)(f^*(x^{\ell_0}_{\langle
i_0\rangle})\neq 0)$.
\smallskip

\noindent We may assume that the group $H$ has universe $\{2i:
i<\lambda\}$. Moreover we may have an increasing continuous sequence
$\lk\gamma_\alpha: \alpha<\lambda\rk$ of limit ordinals such that 
\begin{enumerate}
\item[$\boxtimes_1$] $\{x^k_t: t\in\Dom(I_k)\ \&\ k<k_0\}\subseteq\gamma_0$,
$\alpha_0<\gamma_0$ and $\gamma^k_1<\gamma_0$ if $k_0\leq k<k^*$, $1<n_k$,
and  
\item[$\boxtimes_2$] $H_\alpha\stackrel{\rm def}{=}H\rest \{2i:
i<\gamma_\alpha\}$ is the subgroup of $H$ generated by 
\[\{x^k_t: k<k_0\ \mbox{ or }\ [k_0\leq k<k^*\ \&\ t(0)<\gamma_\alpha]\}.\]
\end{enumerate}
For $k<k^*$ we define the reduction $I_k^{\rm red}$ of the ideal $I_k$ by:
\smallskip

\noindent if $k<k_0$ then $I_k^{\rm red}=I_k$,

\noindent if $k_0\leq k<k^*$ and $n_k>1$ then $I_k^{\rm red}=
I_{\lk\lambda^k_1,\ldots,\lambda^k_{n_k-1}\rk}$ and

\noindent if $k_0\leq k<k^*$ and $n_k=1$ then $I_k^{\rm red}=
I_{\lk\aleph_0\rk}$. 
\smallskip

\noindent Next we define $y^k_t[\gamma_\alpha]$ (for $k<k^*$, $t\in
\Dom(I_k^{\rm red})$) as
\[\begin{array}{ll}
x^k_t &\mbox{ if }\ k<k_0,\\
x^k_{\lk\gamma_\alpha\rk\conc t} &\mbox{ if }\ k_0\leq k<k^*,\
n_k>1,\\ 
x^k_{\lk\gamma_\alpha+t(0)\rk} &\mbox{ if }\ k_0\leq k<k^*, n_k=1.
\end{array}\]
It follows from $\boxtimes_1$, $\boxtimes_2$ and the fact that
$\gamma_\alpha$ are limit ordinals that $y^k_t[\gamma_\alpha]\in 
H_{\alpha+1}$ for all $k<k^*$, $t\in \Dom(I_k^{\rm red})$. The subgroup
generated by these elements with some side elements will be the one to which
we will apply the inductive hypothesis.
\medskip

\noindent Let $E\subseteq \acc(\{\alpha<\lambda: \alpha=\gamma_\alpha\})$ be
a thin enough club of $\lambda$. By our assumptions we find a stationary set
$S\subseteq E$ such that
\begin{enumerate}
\item[$(\alpha)$] for every $\theta=\cf(\theta)<\lambda$, $\{\delta\in S:
\cf(\delta)=\theta\}$ is stationary,
\item[$(\beta)$] $S$ a set of singular limit ordinals,
\item[$(\gamma)$] $\lambda=\mu^+\Rightarrow S\subseteq [\mu+1, \lambda)$ and
if $\lambda$ is inaccessible then $S$ is a set of strong limit singular
cardinals and
\item[$(\delta)$] $S$ does not reflect.
\end{enumerate}
We will use the principle formulated in \ref{principle} to choose by
induction on $\alpha<\lambda$ a group $G_\alpha$ with the universe
$\gamma_\alpha$ and extending $H_\alpha$. For this we have to define finite
vocabularies $\bar{P}$, $\bar{Q}$, $\bar{R}$ and a formula $\varphi$. Thus
we declare that $\bar{P}$ is $\langle P_0,P_1,\ldots\rangle$, $P_0$ a unary
predicate and $P_1$ a unary function symbol, $\bar{Q}=\langle Q_0,\ldots
\rangle$ where $Q_0$ is a binary function symbol, and $\langle M,\bar{Q}^M,
\bar{R}^M\rangle\models\varphi$ means:
\begin{enumerate}
\item[(a)] $\langle M, Q^M_0\rangle$ is a group, $P_0^M$ is its subgroup
(intension: $H$), and $P^M_1\rest P^M_0\in\HOM(P^M_0,\bbz)$ (intension:
$f^*$) 

[we should use some additional predicates to encode $\bbz,\bbz/p\bbz,\ldots$], 
\item[(b)] $\bar R^M$ encodes a homomorphism $f\in\HOM(M,\bbz)$ extending
$f^*$. 
\end{enumerate}
Let functions $\GN^\delta,\GB^\delta_\varepsilon,\Gz$ be given by
\ref{principle} for $\bar{P}$, $\bar{Q}$, $\bar{R}$ and $\varphi$ as defined
above.   
\medskip

Now, by induction on $\alpha<\lambda$, we choose models $M^1_\alpha$ (i.e.\
groups $G_\alpha$, their subgroups $P^{M^1_\alpha}_0=H_\alpha$ and
homomorphisms $f^+_\alpha$, $f^*_\alpha$) and $T_\alpha$, $R_\alpha$,
$h_{(g^*,g^+,C)}$ (for $(g^*,g^+,C)\in R_\alpha$), $R^q_\alpha$ and
$h^+_{(h,C)}$ (for $(h,C)\in R^q_\alpha$) such that:
\begin{enumerate}
\item $\lk G_\alpha:\alpha<\lambda\rk$ is an increasing continuous sequence
of free groups,
\item $P^{M^1_\alpha}_0=H_\alpha\subseteq G_\alpha$, $G_\alpha$ is a (free)
group on $\gamma_\alpha$, 
\item if $\beta<\alpha$ then $G_\alpha/G_\beta$ is free,
\item $G_{\alpha}/H_\alpha$ is free,
\item if $\alpha\notin S$, $\alpha<\beta$ then $G_\beta/(G_\alpha+H_\beta)$
is a free group with a basis of size $\|\gamma_\beta\|$,
\item if $\alpha\in S$, $k'< k^*$, $A\in I_{k'}^{\rm red}$ then
$G_{\alpha+1}/(G_\alpha+H_{\alpha+1}^A)$ is free, where 
$H^A_{\alpha+1}$ is the group generated by all elements $x^k_t$ such
that  $k< k^*$, $t\in\Dom(I_k)$, $x^k_t\in H_{\alpha+1}$ but $(\forall
s\in \Dom(I^{\rm red}_{k'}))(x^k_t= y^{k'}_s[\gamma_\alpha]\
\Rightarrow\ s\in A)$,
\item $f^*_\alpha=f^*\rest H_\alpha$.
\item For $\alpha\in E$ let $N^\alpha=\GB^\alpha(M^1_\alpha)$, $\Gz_\alpha=
\Gz(M^1_\alpha)$, $N^\alpha_i=\GB^\alpha_i(M^1_\alpha)$ for $i<\Gz_\alpha$.

{\em Remark:} Since the group $G_\alpha/H_\alpha$ is free we have
$M^1_\alpha\in\bk^{-}_\alpha$. If $\Gz(M^1_\alpha)=0$ then 
$N^\alpha$ is empty, and below $T_\alpha=R_\alpha=\emptyset$.

Assume $\alpha\in E$ (so $\gamma_\alpha=\alpha$). Then $T_\alpha$ is
the family of all pairs $(g^*,g^+)$ of homomorphisms
$g^*\in\HOM(H_\alpha,\bbz)\cap N^\alpha$,
$g^+\in\HOM(G_{\alpha},\bbz/p\bbz)\cap N^\alpha$ such that
$g^*/p\bbz\subseteq g^+$.
\item $R_\alpha$ is the family of all triples $(g^*,g^+, C)$
such that:
\begin{enumerate}
\item[(a)] $C\in N^\alpha$ is a nonempty closed subset of $\alpha\cap E$,
$(g^*,g^+)\in T_\alpha$, 
\item[(b)] for $\beta\in C$:\ \ \ $(g^*\rest H_\beta,g^+\rest G_\beta,
C\cap\beta)\in R_\beta$, 
\item[(c)] for $\beta<\gamma$ in $C$: $h_{(g^*\rest H_\beta,g^+\rest G_\beta,
C\cap \beta)}\subseteq h_{(g^*\rest H_\gamma, g^+\rest G_\gamma,
C\cap\gamma)}$, 
\item[(d)] if $\beta\in [\min C,\alpha)\cap E$ then for all $k'<k^*$:
\[\{t\in\Dom(I_{k'}^{\rm red}):\ g^*(y^{k'}_t[\gamma_\beta])\neq 0\}\in
I_{k'}^{\rm red}.\]
\end{enumerate}
$R^q_\alpha$ (for $q\in\bP$) is the family of all pairs $(h,C)$
such that:
\begin{enumerate}
\item[(a)] $C\in N^\alpha$ is a nonempty closed subset of $\alpha\cap E$,
$h\in\HOM(H_\alpha,\bbz/q\bbz)\cap N^\alpha$, 
\item[(b)] for $\beta\in C$:\ \ \ $(h\rest H_\beta,C\cap\beta)\in
R^q_\beta$,  
\item[(c)] for $\beta<\gamma$ in $C$: $h^+_{(h\rest H_\beta,C\cap\beta)}
\subseteq h^+_{(h\rest H_\gamma,C\cap\gamma)}$, 
\item[(d)] if $\beta\in [\min C,\alpha)\cap E$ then for all $k'<k^*$:
\[\{t\in\Dom(I_{k'}^{\rm red}):\ h(y^{k'}_t[\gamma_\beta])\neq 0\}\in
I_{k'}^{\rm red}.\]
\end{enumerate}
\item If $(g^*,g^+, C) \in R_\alpha$ {\it then} $h_{(g^*,g^+, C)}\in
\HOM(G_\alpha,\bbz)\cap N^\alpha$, $h_{(g^*,g^+, C)}/p\bbz=g^+$,
$g^*\subseteq h_{(g^*,g^+,C)}$ and 
\[\beta\in C\ \ \Rightarrow\ \ h_{(g^*\rest H_\beta,g^+\rest G_\beta,
C\cap\beta)} \subseteq h_{(g^*,g^+,C)},\]
if $(h,C)\in R^q_\alpha$ then $h^+_{(h,C)}\in\HOM(G_\alpha,\bbz/q\bbz)$
extends $h\cup\bigcup\limits_{\beta\in C} h^+_{h\rest H_\beta,C\cap\beta}$. 
\item If $(g^*,g^+, C)\in R_\alpha$, $g^*_0\in\HOM(H_{\alpha+1},\bbz)$ is
such that $g^*\subseteq g^*_0$ and for every $k'<k^*$:
\[\{t\in\Dom(I_{k'}^{\rm red}):\ g^*_0(y^{k'}_t[\gamma_\alpha])\neq 0\}\in
I_{k'}^{\rm red}\]
and if $g^+_0\in\HOM(G_{\alpha+1},\bbz/p\bbz)$ is such that
$g^*_0/p\bbz\subseteq g^+_0$, $g^+\subseteq g^+_0$

\noindent {\em then} there is $h'\in \HOM(G_{\alpha+1}, \bbz)$ extending
$h_{(g^*,g^+, C)}$ and such that $g^*_0\subseteq h'$, $h'/p\bbz=g^+_0$;

if $(h,C)\in R^q_\alpha$, $h_0\in\HOM(H_{\alpha+1},\bbz/q\bbz)$ is such that
$h\subseteq h_0$ and for every $k'<k^*$
\[\{t\in\Dom(I_{k'}^{\rm red}):\ h_0(y^{k'}_t[\gamma_\alpha])\neq 0\}\in
I_{k'}^{\rm red}\]
{\em then} there is $h'\in \HOM(G_{\alpha+1},\bbz/q\bbz)$ extending
$h^+_{(h, C)}\cup h_0$.
\item Assume that $\alpha\in S$ and
\[\GN(M^1_\alpha)=\langle G_{\alpha},H_\alpha,f^*_\alpha,\ldots,f\rangle,\]
where $f\in\HOM(G_\alpha,\bbz)$.\\
{\em If} there is a free group $G^*$, $G_\alpha\cup H_{\alpha+1}\subseteq
G^*$ such that: $\|G^*\|=\|\gamma_{\alpha+1}\|$, $G^*$ satisfies (2)---(6),
(11) (with $G^*$ playing the role of $G_{\alpha+1}$) and 
\begin{enumerate}
\item[$(*)$] \qquad there is no $g'\in\HOM(G^*,\bbz)$ extending $f\cup
f^*_{\alpha+1}$ 
\end{enumerate}
{\em then} $G_{\alpha+1}$ satisfies $(*)$ too.
\end{enumerate}

The limit stages of the construction are actually determined by the
continuity demands of (1), (7). Concerning the requirements (2)---(5) note
that (2) is preserved because of (3) at previous stages, (3) is preserved
because of (2), (4) is kept due to (5) and the fact that the set $S$ is not
reflecting and finally (5) holds at the limit because of (4) at previous
stages and non reflection of $S$ (see e.g. Proposition~IV.1.7 of 
\cite{EM}). There is some uncertainty in defining $h_{(g^*,g^+,C)}$ for
$(g^*,g^+,C)\in R_\alpha$ (for $\alpha\in E$). However it is possible to
find a suitable $h_{(g^*,g^+,C)}$ since in the most difficult case when
$\sup C<\alpha$, $\sup C\in S$ we may apply first (11) and then
(5). Similarly we handle $h^+_{(h,C)}$. 

\noindent If $\alpha\notin S$ then we choose a group
$G_{\alpha+1}\supseteq H_{\alpha+1}\cup G_\alpha$ such that
$G_{\alpha+1}/(H_{\alpha+1}+G_\alpha)$ is a free group with a basis of size
$\|\gamma_{\alpha+1}\|$.

\noindent If $\alpha\in S$ then condition (12) of the construction describes
$G_{\alpha+1}$. (We will see later that this condition is not empty,
i.e.\ that there is a group $G^*$ as there.) 
\medskip

Thus we have carried out the definition and we may put $G=G_\lambda=
\bigcup\limits_{\alpha<\lambda}G_\alpha$. Let us check that the $G$
satisfies the desired properties (the main point will be the requirement
$(\alpha)$ of the theorem). 

By (2) and (3) the group $G$ is free of cardinality $\lambda$ (and
the set of elements is $\lambda$) and it extends each $G_\alpha$. Due to (4)
the quotient $G/H$ is $\lambda$-free.
\medskip

\noindent{\bf Clause $(\alpha)$ of the assertion}\\
This is a consequence of the condition (12) of the construction. Suppose
that the homomorphism $f^*$ has an extension to a homomorphism in $\HOM(G,
\bbz)$. This means that
\[M^1=\lk G, H, f^*,\ldots\rk\in \bk^-_\lambda.\]
By condition {\bf (i)} of \ref{principle} we find a regular cardinal
$\sigma\leq\lambda$ and 
\[M^2=\lk G, H, f^*,\ldots,f\rk\in \bk^+_\lambda(M^1)\quad\mbox{ and
}\alpha\in S\] 
such that $\GN(M^1_\alpha)=M^2_\alpha\prec M^2$ and
\[\begin{array}{lll}
\cf(\alpha)=\omega_1&\mbox{if}&\lambda=\aleph_2,\\
\cf(\alpha)\neq\cf(\mu)&\mbox{if}&\lambda=\mu^+,\ \mu \mbox{ is an
uncountable limit cardinal,}\\
\cf(\alpha)\notin\{\mu,\cf(\theta)\}&\mbox{if}&\lambda=\mu^+,\
\mu=\theta^+,\ \theta>\aleph_0
  \end{array}\]
(remember ({\bf i})($\alpha$)). Look now at the stage $\alpha$ of the
construction. 

Before we continue with the proof we give a claim which helps us to apply
the inductive hypothesis.

\begin{claim}
\label{claim1}
\begin{enumerate}
\item {\em If} $R\subseteq\HOM(G_\alpha,\bbz)\cup\bigcup\limits_{q\in\bP}
\HOM(G_\alpha,\bbz/q\bbz)$, $2^{\|R\|}<\|\alpha\|$ {\em then} for every
$\beta\in\alpha\setminus S$ large enough there is $x\in G_\alpha$ such that: 
\begin{enumerate}
\item[(a)] $(\forall h\in R)(h(x)=0)$
\item[(b)] $x\in G_{\beta+1}$ is a member of a basis of $G_{\beta+1}$ over
$H_{\beta+1}+G_\beta$. 
\end{enumerate}
\item Suppose that $\mu<\alpha<\mu^+$, $\mu$ is a limit cardinal (so we are
in the case $\lambda=\mu^+$), $R\subseteq\HOM(G_\alpha,\bbz)\cup
\bigcup\limits_{q\in\bP}\HOM(G_\alpha,\bbz/q\bbz)$ and $\|R\|=\mu$. Then for
each $\beta\in(\mu,\alpha)\setminus S$ there exist $x_j\in G_{\beta+1}$ for
$j<\cf(\mu)$ such that:   
\begin{enumerate}
\item[(a)] if $h\in R$ then the set $\{j<\cf(\mu):h(x_j)\neq 0\}$ is bounded
in $\cf(\mu)$,
\item[(b)] $\{x_j:j<\cf(\mu)\}$ can be extended to a basis of $G_{\beta+1}$
over $H_{\beta+1}+G_\beta$.
\end{enumerate}
\item Suppose that $\theta^+=\mu<\alpha<\mu^+$, $R\subseteq\HOM(G_\alpha,
\bbz)\cup\bigcup\limits_{q\in\bP}\HOM(G_\alpha,\bbz/q\bbz)$, $\|R\|=\mu$. 
Then for every $\beta\in(\mu,\alpha)\setminus S$ there exists a sequence
$\lk x_{j,k}: j<\mu, k<\cf(\theta)\rk \subseteq G_{\beta+1}$ such that
\begin{enumerate}
\item[(a)] if $h\in R$ then $\{(j,k)\in\mu\times\cf(\theta): h(x_{j,k})\neq
0\}\in I_{\lk\mu,\cf(\theta)\rk}$,
\item[(b)] $\lk x_{j,k}:j<\mu, k<\cf(\theta)\rk$ can be extended to a basis
of $G_{\beta+1}$ over $H_{\beta+1}+G_\beta$.
\end{enumerate}
\item Suppose $\aleph_1<\alpha<\aleph_2$ (so $\lambda=\aleph_2$),
$R\subseteq\HOM(G_\alpha,\bbz)\cup\bigcup\limits_{q\in\bP}\HOM(G_\alpha,
\bbz/q\bbz)$, $\|R\|=\aleph_0$. Then for each $\beta\in
(\aleph_1,\alpha)\setminus S$ there are $x_\ell\in G_{\beta+1}$ (for $\ell<
\omega$) such that
\begin{enumerate}
\item[(a)] if $h\in R$ then the set  $\{\ell<\omega:h(x_\ell)\neq 0\}$ is
finite,
\item[(b)] $\{x_\ell:\ell<\omega\}$ can be extended to a basis of $G_{\beta+
1}$ over $H_{\beta+1}+G_\beta$.
\end{enumerate}
\end{enumerate}
\end{claim}

\begin{proof}[Proof of the claim]
1)\ \ Let $\beta_0<\alpha$ be such that $\|\gamma_{\beta_0}\|>2^{\|R\|}$
(remember $\alpha\in S\subseteq E$, see the choice of $S$). Let $\beta\in
\alpha\setminus S$, $\beta>\beta_0$. Let $\{y_\xi:\xi<\gamma_{\beta+1}\}$ be
a free basis of $G_{\beta+1}$ over $H_{\beta+1}+G_\beta$ (exists by
condition (5) of the construction). If $R$ is finite then considering first
$\|R\|+1$ elements of the basis we find a respective point $x$ in the group
generated by them. If $R$ is infinite, so $2^{\|R\|}=\|{}^R\bbz\|$, then we
find $\xi_0<\xi_1<\gamma_{\beta+1}$ such that $(\forall h\in
R)(h(y_{\xi_0})=h(y_{\xi_1}))$ and we may put $x=y_{\xi_0}-y_{\xi_1}$.
\medskip

\noindent 2)\ \ We follow exactly the lines of (1), but first we have to
choose an increasing sequence $\lk R_j: j<\cf(\mu)\rk$ such that
$\bigcup\limits_{j<\cf(\mu)} R_j=R$, $\|R_j\|<\|R\|$ (and hence
$2^{\|R_j\|}<\mu$ as $\mu$ is a limit cardinal). Now if $\beta\in
(\mu,\alpha)\setminus S$ then we find $\lk x_j:j<\cf(\mu)\rk\subseteq
G_{\beta+1}$ which can be extended to a basis of $G_{\beta+1}$ over
$H_{\beta+1}+G_\beta$ and such that $(\forall h\in R_j)(h(x_j)=0)$.
\medskip

\noindent 3)\ \ Similarly: first find $\lk R_{j,k}: j<\mu, k<\cf(\theta)\rk$
such that $\|R_{j,k}\|<\theta$, the sequence $\lk\bigcup\limits_{k<\cf(
\theta)}R_{j,k}: j<\mu\rk$ is increasing, for each $j<\mu$ the sequence $\lk
R_{j,k}: k<\cf(\theta)\rk$ is increasing and $\bigcup\limits_{j<\mu}
\bigcup\limits_{k<\cf(\theta)}R_{j,k}= R$. Next follow as in 2). 
\medskip

\noindent 4)\ \ Represent $R$ as an increasing (countable) union of finite
sets and follow as in 2) above. The claim is proved. 
\end{proof}

Now we are going back to the proof of clause $(\alpha)$. The following claim
will finish it.

\begin{claim}
\label{claim2}
Suppose that $\alpha$, $f\ldots$ are as chosen earlier. Then there exists a
free group $G^*\supseteq H_{\alpha+1}\cup G_\alpha$ such that $\|G^*\|=\|
\gamma_{\alpha+1}\|$, $G^*$ satisfies the conditions (2)--(6), (11) of the
construction as $G_{\alpha+1}$ and there is no $g'\in\HOM(G^*,\bbz)$
extending $f\rest G_\alpha\cup f^*_{\alpha+1}$. 
\end{claim}

\begin{proof}[Proof of the claim]
Let $R=\{h_{(g^*,g^+,C)}: (g^*,g^+,C)\in R_\alpha\}\cup\{h^+_{(h,C)}:
(h,C)\in R^q_\alpha,\ q\in\bP\}$. By clauses {\bf (g)} and {\bf (b)} of
\ref{principle} we have $\|R\|\leq\|\alpha\|$ (of course $R$ may be empty). 
Let $\lk \alpha_i: i<\cf(\alpha)\rk$ be an increasing continuous sequence
cofinal in $\alpha$ and disjoint from $S$ (possible by the choice of $S$).
\medskip

\begin{description}
\item[\bf Case A] $\alpha$ is a strongly limit singular cardinal (so we are
in the case when $\lambda$ is inaccessible).
\end{description}
We find an increasing sequence $\lk R_i^*: i<\cf(\alpha)\rk$ such that
$\bigcup\limits_{i<\cf(\alpha)} R^*_i=R$ and $\|R^*_i\|<\|\alpha\|$. But in
this case we have 
\[(\forall i<\cf(\alpha))(2^{\|R^*_i\|}<\alpha).\]
So we may apply claim \ref{claim1}(1) to choose by induction on
$i<\cf(\alpha)$ an increasing sequence $\lk j_i: i<\cf(\alpha)\rk\subseteq
\cf(\alpha)$ and $x^{k^*}_i\in G_{\alpha_{j_i}+1}$ such that:
\begin{enumerate}
\item[(a)] $h\in R^*_i\Rightarrow h(x^{k^*}_i)=0$
\item[(b)] $x^{k^*}_i$ is a member of a basis of $G_{\alpha_{j_i}+1}$ over
$H_{\alpha_{j_i}+1}+G_{\alpha_{j_i}}$. 
\end{enumerate}
Since $\lk \alpha_i: i<\cf(\alpha)\rk \subseteq\alpha\setminus S$ is
increasing continuous (and cofinal in $\alpha$) we get that $\{x^{k^*}_i:
i<\cf(\alpha)\}$ can be extended to a basis of $G_\alpha$ over $H_\alpha$.
Now we apply the inductive hypothesis to $k^*+1$, $I^{\rm red}_k$ (for
$k<k^*$), $I_{\lk\cf(\alpha)\rk}$, the group $H^*$ generated by
\[\{y^k_t[\gamma_\alpha]: k<k^*, t\in\Dom(I^{\rm red}_k)\}\cup \{x^{k^*}_i:
i<\cf(\alpha)\}\] 
and the function $(f\cup f^*_{\alpha+1})\rest H^*$. This gives us a group
$G^*_0\supseteq H^*$. Let $H'$ be such that $G_\alpha+H_{\alpha+1}= H^*
\oplus H'$. Then put $G^*=G^*_0\oplus H'$. It satisfies the requirements of
the claim: condition (3) follows from the presence of the $y^k_t[
\gamma_\alpha]$'s part of $H^*$ (remember the inductive assumption
\ref{1.2B}($\beta$)), condition (4) holds due to the $x^{k^*}_i$. It follows
from the fact that the $\alpha_{j_i}$ are cofinal in $\alpha$ (and from the
choice of $x^{k^*}_i\in G_{\alpha_{j_i}+1}$) that (5) is satisfied. 
Similarly, (11) is a consequence of the choice of $x^{k^*}_i$ and the
inductive hypothesis \ref{1.2B}($\gamma$, $\delta$). Finally clause (6)
follows from the inductive assumption \ref{1.2B}($\beta$)
\medskip

\begin{description}
\item[\bf Case B] $\aleph_0<\alpha<\aleph_1$ (so $\lambda=\aleph_1$).
\end{description}
Thus $R$ is at most countable, so let $R=\bigcup\limits_{\ell<\omega}
R_\ell$, where $R_\ell$ are finite increasing with $\ell$. Apply
\ref{claim1}(1) to find an increasing sequence $\langle j_\ell:\ell<\omega
\rangle\subseteq\omega$ and $x^{k^*}_\ell\in G_{\alpha_{j_\ell}+1}$ such
that 
\begin{enumerate}
\item[(a)] $h\in R_\ell\quad\Rightarrow\quad h(x^{k^*}_\ell)=0$,
\item[(b)] $x^{k^*}_\ell$ is a member of a basis of $G_{\alpha_{j_\ell}+1}$
over $H_{\alpha_{j_\ell}+1}+G_{\alpha_{j_\ell}}$.
\end{enumerate}
Proceed as in Case A (so apply the inductive hypothesis to $I^{\rm red}_k$
(for $k<k^*$) and $I_{\langle\aleph_0\rangle}$). 
\medskip

\begin{description}
\item[\bf Case C] $\alpha\in (\mu,\mu^+)$ for some limit cardinal
$\mu>\aleph_0$ (so $\lambda=\mu^+$).
\end{description}
Then we have $\|R\|\leq\mu$ and by claim \ref{claim1}(2) we can choose
$x^{k^*}_{i,j}\in G_{\alpha_i+1}$ (for $i<\cf(\alpha)$, $j<\cf(\mu)$) such
that: 
\begin{enumerate}
\item[(a)] for each $h\in R$ for every $i<\cf(\alpha)$ the set
$\{j<\cf(\mu): h(x^{k^*}_{i,j})\neq 0\}$ is bounded, 
\item[(b)] for each $i<\cf(\alpha)$ the set
$\{x^{k^*}_{i,j}:j< \cf(\mu)\}$ extends to a basis of $G_{\alpha_i+1}$ over
$H_{\alpha_i+1}+G_{\alpha_i}$. 
\end{enumerate}
Now apply the inductive hypothesis for $k^*+1$, $I^{\rm red}_k$ (for
$k<k^*$) and $I_{\lk\cf(\alpha),\cf(\mu)\rk}$ (remember that $\cf(\alpha)
\neq\cf(\mu)$ in this case). 
\medskip

\begin{description}
\item[\bf Case D] $\aleph_1<\alpha<\aleph_2$ (so $\lambda=\aleph_2$ and
$\cf(\alpha)=\omega_1$). 
\end{description}
Write $R$ as an increasing union $\bigcup\limits_{i<\omega_1}R_i$ of
countable sets. Using \ref{claim1}(4) choose $x^{k^*}_{i,\ell}$ (for
$i<\omega_1$, $\ell<\omega$) such that
\begin{enumerate}
\item[(a)] for each $h\in R$, for every sufficiently large $i<\omega_1$ the
set $\{\ell<\omega: h(x^{k^*}_{i,\ell})\neq 0\}$ is finite,
\item[(b)] for each $i<\omega_1$ the set $\{x^{k^*}_{i,\ell}:\ell<\omega\}$
can be extended to a basis of $G_{\alpha_i+1}$ over $H_{\alpha_i+1}+
G_{\alpha_i}$.
\end{enumerate}
Proceed as above (using $I^{\rm red}_k$ (for $k<k^*$) and $I_{\langle
\aleph_1,\aleph_0\rangle}$). 
\medskip

\begin{description}
\item[\bf Case E]  $\alpha\in (\mu,\mu^+)$ for some cardinal $\mu$
such that $\mu=\theta^+>\aleph_1$ (so $\lambda=\mu^+$).
\end{description}
Using claim \ref{claim1}(3) we choose a sequence $\lk x^{k^*}_{l,j,i}:
l<\cf(\alpha), j<\mu, i<\cf(\theta)\rk$ such that 
\begin{description}
\item[(a)] for each $h\in R$ and for every $l<\cf(\alpha)$, for every
$j<\mu$ large enough for every $i<\cf(\theta)$ large enough,
$h(x^{k^*}_{l,j,i})=0$. 
\item[(b)] for every $l<\cf(\alpha)$ the set
$\{x^{k^*}_{l,j,i}: j<\mu,i<\cf(\theta)\}$ can be extended to a basis of
$G_{\alpha_l+1}$ over $H_{\alpha_l+1}+G_{\alpha_l}$. 
\end{description}
Now apply the inductive hypothesis to $k^*+1$, $I^{\rm red}_k$ ($k<k^*$) and
$I_{\lk\cf(\alpha),\mu,\cf(\theta)\rk}$ (remember $\cf(\alpha)\notin\{\mu,
\cf(\theta)\}$ in this case). 
\medskip

\noindent This completes the proof of claim \ref{claim2}. 
\end{proof}

It follows from the above claim that at the stage $\alpha$ of the
construction we had a nontrivial application of the condition (12)
``killing'' the function $f$. This gives a contradiction proving the clause
$(\alpha)$. 
\medskip

\noindent{\bf Clause $(\beta)$ of the assertion}\\
It follows from conditions (6) and (5) of the construction.
\medskip

\noindent{\bf Clause $(\gamma)$ of the assertion}\\
Assume that $g^*$, $g^+$ are as there. Then by the clause {\bf (j)} of
\ref{principle} we have a club $C\subseteq E$ such that for each $\alpha\in C$ 
\[(g^*\rest H_\alpha,g^+\rest G_\alpha,C\cap\alpha) \in R_\alpha.\]
Consequently we may use the functions $h_{(g^*\rest H_\alpha,g^+\rest
G_\alpha,C\cap\alpha)}$ for $\alpha\in C$.
\medskip

\noindent{\bf Clause $(\delta)$ of the assertion}\\
Like clause $(\gamma)$.
\end{proof}

Before we state the main result let us recall basic properties of $\Ext$.
First note that

\begin{quotation}
\noindent if $G$ is an (abelian) group satisfying
$G\models(\forall x)(px=0)$ 

\noindent then $G$ is a vector space over $\bbz/p\bbz$.
\end{quotation}

\begin{definition}
\begin{enumerate}
\item For a group $G$ and $p\in \bP$ let $\nu_p(G)$ be the dimension of
$\Ext_p(G, \bbz)$ as a vector space over $\bbz/p\bbz$ where 
\[\Ext_p(G, \bbz)=\{x\in \Ext(G, \bbz): \Ext(G,\bbz)\models px=0\}.\]
\item For a group $G$ and let $\nu_0(G)$ be the rank (=maximal cardinality of
an independent subset) of the torsion free group $\Ext(G,
\bbz)/tor(\Ext(G,\bbz))$ where for a group $G'$:
\[tor(G')=\{x\in G':\mbox{ for some } n,\ 0<n\in\bbz\mbox{ we have }G'\models
n x=0\}.\]
\end{enumerate}
\end{definition}

\begin{lemma}
[see Fuchs \cite{Fu} or Eklof Mekler {\cite[Ch.\ XII]{EM}}] \ \\
Let $G$ be an abelian torsion-free group. Then:
\begin{enumerate}
\item For $p\in \bP$, $\nu_p(G)$ is the dimension of the vector space
\[\HOM(G, \bbz/p\bbz)/\HOM^-(G, \bbz/p\bbz)\]
over the field $\bbz/p\bbz$,
where $\HOM^-(G, \bbz/p\bbz)\stackrel{\rm def}{=}\{f/p:f\in \HOM(G, \bbz)\}$.
\item $\Ext(G, \bbz)$ is a divisible group, hence characterized up to
isomorphism by cardinals $\nu_0(G), \nu_p(G)$ (for $p\in\bP$).
\end{enumerate}
\end{lemma}

\begin{theorem}
[Hiller, Huber, Shelah \cite{HHSh91}] {\em ({\bf V}={\bf L})}

\noindent If a group $G$ is not free, moreover it is not $G_1\oplus G_2$
with $G_2$ free, $\| G_1\| <\| G\|$  {\em then} $\nu_0(G)=2^{\|G\|}$. 
\end{theorem}

\begin{remark}
If $G=G_1\oplus G_2$ and $G_2$ is
free then $\Ext(G,\bbz)\cong\Ext(G_1,\bbz)$ so the demand is reasonable.
\end{remark}

\begin{maintheorem}[V=L]
\label{1.2}
Suppose that $\lambda$ is an uncountable regular cardinal which is smaller
then the first weakly compact cardinal. Let $\lambda_p\le \lambda^+$ for
$p\in \bP$. {\em Then} there exists a (torsion free) strongly $\lambda$-free
group $G$ such that $\|G\|=\lambda$, $\nu_p(G)=\lambda_p$, and $\nu_0(G)=
\lambda^+$.
\end{maintheorem}

\begin{proof}
During the proof we will use consequences of the assumption ${\bf V}={\bf
L}$ like GCH, the principle proved in \ref{principle} etc without recalling
the main assumption. 

The construction is much easier if $\lambda_p=\lambda^+$ for some $p\in\bP$
and $\lambda_q=0$ for all $q\neq p$ (remember that $\Ext(\bigoplus\limits_{n
\in\omega} G_n,\bbz)=\prod\limits_{n\in\omega}\Ext(G_n,\bbz)$). Therefore we
assume that we are done with this particular case and we assume that
$\lambda_p\leq\lambda$ for all $p\in\bP$. 

We shall build a $\lambda$-free group $G=G_\lambda=\bigcup\limits_{\alpha<
\lambda}G_\alpha$ with universe $\lambda$ (the sequence $\langle G_\alpha:
\alpha<\lambda\rangle$ increasing continuous, $G_\alpha$ a group on an
ordinal $\gamma_\alpha<\lambda$ for $\alpha <\lambda$). As witnesses for
$\nu_p(G)\ge \lambda_p$  there will be also  homomorphisms $f^{p,
\zeta}_\lambda\in\HOM(G,\bbz/p\bbz)$ for $\zeta<\lambda_p$, $f^{p,
\zeta}_\lambda=\bigcup_{\alpha<\lambda} f^{p,\zeta}_\alpha$. For the
witnesses to work we need: 
\begin{enumerate}
\item[$(*)_1$] {\em If\/} $p\in \bP$, $0<n<\omega$, $\zeta^0<\dots<
\zeta^{n-1}<\lambda_p$, $a_\ell\in\{1/p\bbz,\dots,(p-1)/p\bbz\}$ (for
$\ell<n)$,\\
{\em then\/} $\sum_{\ell<n} a_\ell f^{p,\zeta^\ell}_\lambda\notin\HOM^-(G,
\bbz/p\bbz)$,\\
(of course $\sum_{\ell<n} a_\ell f^{p,\zeta^\ell}_\lambda\in\HOM(G,\bbz/
p\bbz))$.
\end{enumerate}
This is equivalent to
\begin{enumerate}
\item[$(*)_2$] there are no $p\in \bP$, $0<n<\omega$, $\zeta^0<\dots<
\zeta^{n-1}<\lambda_p$, $a_\ell\in \{1/p\bbz, \dots, (p-1)/p\bbz\}$ (for
$\ell<n$) and $g\in\HOM(G,\bbz)$ such that $g/p\bbz=\sum_{\ell<n} a_\ell
f^{p, \zeta^\ell}_\lambda$.
\end{enumerate}
We shall also have to take care showing that $\nu_p(G)$ is not $>\lambda_p$
(if $\lambda_p<2^\lambda$) and for this it suffices to show that
$\{f^{p,\zeta}_\lambda:\zeta<\lambda_p\}$ generates $\HOM(G,\bbz/p\bbz)$ over
$\HOM^-(G, \bbz/p\bbz)$. For this we shall use the $h_{(g, C)}$ (and
$T^p_\alpha$) below.

By induction on $\alpha<\lambda$ choose an increasing continuous sequence
$\langle\gamma_\alpha:\alpha<\lambda\rangle\subseteq\lambda$ such that
$\lambda_p<\lambda\ \ \Rightarrow\ \ \lambda_p+\omega<\gamma_0$,
$\gamma_{\alpha+1}=\gamma_\alpha+\gamma_\alpha$.  

For our given $\lambda$, we want to use \ref{principle}; we use a club
$E\subseteq\acc(\{\alpha<\lambda:\gamma_\alpha=\alpha\})$ thin enough. As
${\bf V}={\bf L}$ we find a stationary set $S\subseteq E$ such that:
\begin{enumerate}
\item[$(\alpha)$] for every $\theta=\cf(\theta)<\lambda$, $\{\delta\in S:
\cf(\delta)=\theta\}$ is stationary,
\item[$(\beta)$] $S$ a set of singular limit ordinals,
\item[$(\gamma)$] $\lambda=\mu^+\Rightarrow S\subseteq [\mu+1, \lambda)$ and
if $\lambda$ is inaccessible then $S$ is a set of strong limit singular
cardinals and
\item[$(\delta)$] $S$ does not reflect.
\end{enumerate}
Now, $\bar P$ is empty, $\bar{Q}=\langle Q_0,Q_1,Q_2,Q_3,Q_4\ldots\rangle$
where $Q_0$ is a binary function symbol, $Q_1$, $Q_4$ are 3-place ones and
$Q_2$, $Q_3$ are binary predicates and $\langle M, \bar Q^M, \bar
R^M\rangle\models\varphi$ means: 
\begin{enumerate}
\item[(a)] $\langle M, Q^M_0\rangle$ is a group, $Q^M_1(p,\zeta,\cdot)$ is a
homomorphism from the group to ${\mathbb Z}/p{\mathbb Z}$ with $p,\zeta$
variable (it corresponds to $f^{p,\zeta}$); also $\bbz,\bP,\ldots$ are coded
in some way (see below),
\item[(b)] $\bar R^M$ codes a counterexample to $(*)_1$, i.e. $p$, $n$,
$\zeta^0,\ldots,\zeta^{n-1}$, $a_0,\ldots,a_{n-1}$, $f$, $g$ such that
$g\in\HOM(M,\bbz)$, $f=\sum\limits_{\ell<n}a_\ell f^{p,\zeta^\ell}=g/p\bbz
\in\HOM(M,\bbz/p\bbz)$.
\end{enumerate}
Let $\GN^\delta, \GB^\delta_\varepsilon, \Gz$ be as gotten in \ref{principle}. 
Choose a sequence $\langle A_\alpha:\alpha<\lambda\rangle$ such that
$A_\alpha\subseteq\alpha$ for $\alpha<\lambda$ and 
\smallskip

if $A\subseteq\beta$, $\beta<\lambda$

then there is $\alpha\in (\beta, (\|\beta\|+\aleph_0)^+)$ such that $A=
A_\alpha$ 
\smallskip

\noindent (remember we have GCH).

We now choose by induction on $\alpha<\lambda$ the following objects:\quad
$M^1_\alpha$ (i.e. a group $G_\alpha$ and homomorphisms $f^{p,\zeta}_\alpha$
(for $p\in \bP, \zeta\in \lambda_p\cap \alpha)$ and $Q_2^{M^1_\alpha}$,
$Q_3^{M^1_\alpha}$, $Q_4^{M^1_\alpha},\ldots$), $T^p_\alpha$, $R^p_\alpha$,
$h_{(g,C)}$ (for $(g,C)\in R^p_\alpha$) such that:
\begin{enumerate}
\item $G_\alpha$ is a free group with universe $\gamma_\alpha$,
\item $G_\alpha$ is increasing continuous in $\alpha$,
\item if $\beta<\alpha, \beta\notin S$ then $G_\alpha/G_\beta$ is a free
group of size $\|\gamma_\alpha\|$,
\item $f^{p, \zeta}_\alpha\in \HOM(G_\alpha, \bbz/p\bbz)$, $M^1_\alpha$ is
$(G_\alpha, f^{p, \zeta}_\alpha)$ considering $f^{p, \zeta}_\alpha(x)$ a
function with three places (so $f^{p,\zeta}_\alpha$ is not defined for
$\zeta\geq\gamma_\alpha$),
\item if $\beta<\alpha$ then $f^{p, \zeta}_\beta\subseteq f^{p,
\zeta}_\alpha$ (so that $f^{p, \zeta}_\alpha$ is increasing continuous in
$\alpha$), 
\item if $\alpha\notin S$ then there is a basis $Y_\alpha$ of $G_{\alpha+1}$
over $G_\alpha$ such that: 

{\em if} $p\in \bP$, $n<\omega$, $\zeta^0<\dots<\zeta^{n-1}<\lambda_p\cap
\gamma_\alpha$ and $a_0,\ldots,a_{n-1}\in \bbz/p\bbz$ 

{\em then} for $\|\gamma_\alpha\|$ members $y\in Y_\alpha$, $f^{p,
\zeta^\ell}_{\alpha+1}(y)=a_\ell$ (for $\ell<n$) and $f^{q,\zeta}_{\alpha+1}
(y)=0$ if $q\in\bP$, $\zeta<\lambda_q\cap\gamma_\alpha$, $(q,\zeta)\notin\{
(p,\zeta^0),\ldots,(p,\zeta^{n-1})\}$. 
\item Let $\alpha\in E$, $N^\alpha=\GB^\alpha(M^1_\alpha)$, $\Gz_\alpha=
\Gz(M^1_\alpha)$, $N^\alpha_i=\GB^\alpha_i(M^1_\alpha)$ for $i<\Gz_\alpha$.

Note: $M^1_\alpha\in\bk^{-}_\alpha$ by clause (1), the universe of
$N^\alpha$ is a transitive set (see {\bf (f)} of \ref{principle}) so 
$\ord\cap N^\alpha$ is an ordinal greater then $\alpha$ (if non zero).

Assume $\alpha\in S$ (so $\gamma_\alpha=\alpha)$ and $G_\alpha$, $\lk f^{p,
\zeta}_\alpha:\zeta<\lambda_p\cap \alpha\rk$ belong to $N^\alpha$. Then we
choose by induction on $\zeta\in \ord\cap N^\alpha\setminus(\lambda_p\cap
\alpha)$ the function $f^{p, \zeta}_\alpha\in N^\alpha$ as the
$<^*_{N^\alpha}$-first member of $\HOM(G_\alpha,\bbz/p\bbz)$ (as a vector
space over the field $\bbz/p\bbz$) which does not depend on
$\{f^{p,\xi}_\alpha:\xi<\zeta\}$. 

Let $f^{p, \zeta}_\alpha$ be defined if and only if $\zeta<\zeta(p,\alpha)$. 

Note: if $\Gz(M^1_\alpha)=0$ then $\zeta(p,\alpha)=\lambda_p\cap\alpha$,
$N^\alpha$ is empty, and below $T^p_\alpha=R^p_\alpha=\emptyset$.
\item  $T^p_\alpha=\{f^{p, \xi}_\alpha:\lambda_p\cap \alpha\le \xi<\zeta(p,
\alpha)\}$.
\item $R^p_\alpha$ is the family of all pairs $(f^{p,\xi}_\alpha, C)$
such that:
\begin{enumerate}
\item[(a)] $f^{p,\xi}_\alpha\in T^p_\alpha$ and $C\in N^\alpha$ is a closed
subset of $\alpha\cap E$, 
\item[(b)] for $\beta\in C$:\ \ \ $f^{p, \xi}_\alpha\rest G_\beta\in
T^p_\beta$ and $C\cap \beta\in N^\beta$,
\item[(c)] for $\beta<\gamma$ in $C$, $h_{(f^{p, \xi}_\alpha\rest G_\beta,
C\cap \beta)}\subseteq h_{(f^{p, \xi}_\alpha\rest G_\gamma, C\cap\gamma)}$,
\end{enumerate}
\item if $(g,C) \in R^p_\alpha$ {\it then} $h_{(g, C)}\in \HOM(G_\alpha,
\bbz)\cap N^\alpha$ and $h_{(g, C)}/p\bbz=g$ and $\bigcup\limits_{\beta\in
C} h_{(g\rest G_\beta,C\cap\beta)}\subseteq h_{(g,C)}$,
\item if $(g, C)\in R^p_\alpha$, $g\subseteq g'\in\HOM(G_{\alpha+1},
\bbz/p\bbz)$\\
{\em then} there is $h'\in \HOM(G_{\alpha+1}, \bbz)$ extending
$h_{(g, C)}$ and $h'/p\bbz=g'$.
\item Assume that $\alpha\in S$, $0<n<\omega$ and
\[\GN(M^1_\alpha)=\langle G_{\alpha},f^{p,\zeta}_\alpha,\ldots,p,n,\zeta^0,
\ldots,\zeta^{n-1},a_0,\ldots,a_{n-1},f,g\rangle,\]
where $f=\sum_{\ell<n} a_\ell f^{p,\zeta^\ell}_\alpha\in\HOM(G_\alpha,\bbz/
p\bbz)$, $a_\ell\in \bbz/p\bbz\setminus\{0\}$, $f=g/p\bbz$,
$g\in\HOM(G_\alpha,\bbz)$.

{\em If\/} there is a free group $H$, $G_\alpha\subseteq H$ such that:
for every $q\in\bP$, $\zeta<\lambda_q\cap\gamma_\alpha$ the homomorphism
$f^{q,\zeta}_\alpha$ can be extended to a member of $\HOM(H,\bbz/q\bbz)$,
and the quotient $H/G_\beta$ is free for $\beta\in(\alpha\setminus S)$,
$\|H\|=\|\gamma_{\alpha+1}\|$, and $H$ satisfies (11) (with $H$ playing the
role of $G_{\alpha+1}$) and there are some $f^{p,\zeta}_{*}$ satisfying
(4),(5) (with $H$ as $G_{\alpha+1}$ and $f^{p,\zeta}_{*}$ as $f^{p,\zeta}_{
\alpha+1}$) such that for $f'=\sum\limits_{\ell<n}a_\ell f^{p,\zeta}_*\in
\HOM(H,\bbz/p\bbz)$ we have:
\begin{enumerate}
\item[$(*)$] $\quad\neg (\exists g')[g\subseteq g'\in\HOM(H,\bbz)\ \ \&\
\ g'/p\bbz=f']$
\end{enumerate}

{\em then\/} $H=G_{\alpha+1}$, $f'=\sum_{\ell<n} a_\ell f^{p,
\zeta^\ell}_{\alpha+1}$ satisfy $(*)$ too.
\item {\em if\/} $\bigwedge\limits_{p\in\bP}\lambda_p=0$, (so (12) is an
empty demand), $\alpha\in S$, and there is a group $H$ such that $G_\alpha
\subseteq H$ and for each $\beta\in \alpha\setminus S$ the quotient
$H/G_\beta$ free, and $\|H\|=\|\gamma_{\alpha+1}\|$, and it satisfies (11)
(with $H$ playing the role of $G_{\alpha+1}$) and $H/G_\alpha$ is not free

{\em then\/} $G_{\alpha+1}/G_\alpha$ is not free.

\item $Q_2^{M^1_\alpha},Q_3^{M^1_\alpha}\subseteq\gamma_\alpha\times
\gamma_\alpha$ are such that for each $\beta<\gamma_\alpha$ we have:
\begin{quotation}
\noindent $A_\beta=\{i<\gamma_\alpha: M^1_\alpha\models Q_2(\beta,i)\}$ and

\noindent if $\beta$ is limit then $\{i<\gamma_\alpha: M^1_\alpha\models
Q_3(\beta,i)\}$ is a cofinal subset of $\beta$ of the order type $\cf(\beta)$.
\end{quotation}
$Q_4^{M^1_\alpha}$ is such that if $\zeta,\xi<\gamma_\alpha$, $\|\zeta\|=
\|\xi\|$ then the function 
\[Q_4^{M^1_\alpha}(\zeta,\xi,\cdot)\rest\zeta:\zeta\longrightarrow\xi\]
is one-to-one and onto. We require that $Q^{M^1_\alpha}_2,Q^{M^1_\alpha}_3$
and $Q^{M^1_\alpha}_4$ are increasing with $\alpha$, of course.
\end{enumerate}

The conditions (6), (12), (13) and (14) fully describe what happens at
successor stages of the construction. Limit cases are determined by the
continuity demands (2) and (5). Note that the demands (1), (3) are preserved
at the limit stages as the set $S$ is not reflecting (see e.g.\ 
\cite[Proposition~IV.1.7]{EM}). Hence there is no problem to carry out the
definition and let $G=G_\lambda=\bigcup_{\alpha<\lambda}G_\alpha$ (though it
is not so immediate that $G$ is not free!).
\medskip

Some of the desired properties are clear:
\medskip
\begin{enumerate}
\item[$(\otimes)_1$]  $G$ is a group of cardinality $\lambda$ (and the set of
elements is $\lambda$) extending each $G_\alpha$
\end{enumerate}
(by (1)+(2)),

\begin{enumerate}
\item[$(\otimes)_2$] $G$ is $\lambda$-free and even strongly $\lambda$-free
\end{enumerate}
(by (1) each $G_\alpha$ is free so $G$ is $\lambda$-free; by (3) if $\beta
\in\lambda\setminus S$, $\beta<\alpha<\lambda$ then $G_\alpha/G_\beta$ is
free, so $G$ is strongly $\lambda$-free, see e.g. \cite[pp 87--88]{EM}).

Let $f^{p,\zeta}=f^{p,\zeta}_\lambda=\bigcup_{\alpha<\lambda}f_\alpha^{p,
\zeta}$ for $p\in \bP$, $\zeta<\lambda_p$. 

\begin{enumerate}
\item[$(\otimes)_3$] $f^{p, \zeta}\in \HOM(G, \bbz/p\bbz)$ (extending each
$f^{p, \zeta}_\alpha$ for $\alpha<\lambda$) 
\end{enumerate}
(by (4)+(5)). Before checking the main properties of $G$ let us note the
following two facts (which explain the condition (14) of the construction).
\begin{claim}
\label{cl1}
If $\delta\in E$, $\Gz(M^1_\delta)>0$ and $\kappa<\delta$, $\|\kappa\|<\|
\delta\|$ then ${\mathcal P}(\kappa)\subseteq N^\delta$.
\end{claim}

\begin{proof}[Proof of the claim] By the clauses {\bf (f)}, {\bf (h)} of
\ref{principle} and condition (14) of the construction we have $A_i\in
N^\delta_0$ for all $i\in(\kappa, \|\kappa\|^+)$. By the choice of the
sequence $\langle A_i: i<\lambda\rangle$ we are done. 
\end{proof}

\begin{claim}
\label{cl1A}
If $\delta\in E$, $\Gz(M^1_\delta)>0$ then there is an increasing
cofinal in $\delta$ sequence $\langle \beta^\delta_i: i<\cf(\delta)\rangle$
such that for every $i^*<\cf(\delta)$ we have 
\[\langle\beta^\delta_i: i<i^*\rangle\in N^\delta.\]
\end{claim}

\begin{proof}[Proof of the claim]
By the clauses {\bf (k)}, {\bf (l)} of \ref{principle} we have a club
$e_\delta\subseteq\delta$ such that $\otp(e_\delta)<\delta$ and for each
$\alpha<\delta$ the intersection $e_\delta\cap \alpha$ is in $N^\delta$. The
set $b^*\stackrel{\rm def}{=} \{i<\gamma_\delta: M^1_\delta\models
Q_3(\otp(e_\delta), i)\}$ is an increasing cofinal subset of
$\otp(e_\delta)$ of the order type $\cf(\delta)=\cf(\otp(e_\delta))$. It
follows from the condition {\bf (h)} of \ref{principle} that $b^*\in
N^\delta$. But with $e_\delta$ and $b^*$ in hands we may easily build
$\langle\beta^\delta_i: i<\cf(\delta)\rangle$ as required.
\end{proof}

Now comes the main point:

\begin{enumerate}
\item[$(\otimes)_4$]  if $p\in \bP$, $0<n<\omega$, $\zeta^0<\dots<\zeta^{n-1}<
\lambda_p$, $a_0, \dots, a_{n-1}\in \bbz/p\bbz\setminus \{0\}$,
$f=\sum_{\ell<n}a_\ell f^{p, \zeta^\ell}\in \HOM(G, \bbz/p\bbz)$\\
{\em then} $f\notin \HOM^-(G, \bbz/p\bbz)$.
\end{enumerate}
\bigskip

\centerline{\em Why $(\otimes)_4$?}
\bigskip

\noindent Assume that $(\otimes)_4$ fails, so there are $p\in\bP$, $\zeta^0<
\ldots<\zeta^{n-1}$, $a_0,\ldots,a_{n-1}\in\bbz/p\bbz\setminus\{0\}$ and
$g\in \HOM (G,\bbz)$ with $f=\sum\limits_{\ell<n}a_\ell f^{p,\zeta^\ell}=
g/p\bbz$. Let
\[M^2=\langle G,f^{p,\zeta},\ldots,p,n,\zeta^0,\ldots,\zeta^{n-1},
a_0,\ldots,a_{n-1},f,g\rangle, \quad M^2_\delta=M^2\rest\delta.\]
By \ref{principle} condition ({\bf i}) without loss of generality
(i.e.\ possibly replacing
\[p,n,\zeta^0,\ldots,\zeta^{n-1},a_0,\ldots,a_{n-1},f,g\]
by some other
\[p^*,n^*,\zeta^0_*,\ldots,\zeta^{n^*-1}_*,a_0^*,\ldots,a_{n^*-1}^*,f^*,g^*\]
with the same properties) we have: \quad the set
\[\begin{array}{rr}
S^*\stackrel{\rm def}{=}\{\delta\in S:&\GN(M^1_\delta)=M^2_\delta \mbox{ and
} M^2_\delta\prec M^2\mbox{ and }\ \\
 &\Gz(M^1_\delta)=\cf(\delta) \mbox{ or } \Gz(M^1_\delta)=0\ \}
  \end{array}\]
is stationary. (Just applying {\bf (i)} choose $\theta=\sigma$ when
$0<\sigma<\lambda$ and that arbitrary regular $\theta<\lambda$ in other
cases.) Choose $\delta\in S^*$. Remember that $S^*\subseteq S$, so
e.g. $\delta=\gamma_\delta$, $\cf(\delta)<\delta$ and  
\smallskip

$\lambda=\mu^+\Rightarrow\delta\in [\mu+1, \lambda)$ and

if $\lambda$ is inaccessible then $\delta$ is a strongly limit singular
cardinal. 
\smallskip

\noindent Let us first consider the case $\Gz(M^1_\delta)\neq 0$ (so
$\cf(\delta)=\Gz(M^1_\delta)$). To show $(\otimes_4)$ we will need the
following technical but useful claims.

\begin{claim}
\label{cl2}
{\em If} $R\in N^\delta$, $R\subseteq\HOM(G_\delta,\bbz)\cup\bigcup\limits_{
q\in\bP}\HOM(G_\delta,\bbz/q\bbz)$, $2^{\|R\|}<\|\delta\|$ (so
$2^{2^{\|R\|}}\leq\delta$) {\em then} for every $\beta\in \delta\setminus S$
large enough there is $x\in G_\delta$ such that:
\begin{enumerate}
\item[(a)] $h\in R\Rightarrow h(x)=0$,
\item[(b)] $x\in G_{\beta+1}$, moreover $G_\beta\oplus(\bbz x)$ is a direct
summand of $G_{\beta+1}$,
\item[(c)] $g(x)\neq  0$.
\end{enumerate}
\end{claim}

\begin{proof}[Proof of the claim]
First assume that $R$ is infinite, so $2^{\|R\|}=\|{}^R \bbz\|$.

\noindent As $2^{\|R\|}<\|\delta\|$ clearly there are $x\in
G_\delta\setminus\{0\}$ satisfying (a). If $x$ satisfies (a)+(c), $x\in
G_\beta$ and $\beta\in \delta\setminus S$ large enough then we can find $y\in
G_{\beta+1}$ which is a member of a basis of $G_{\beta+1}$ over $G_\beta$ and
which satisfies (a) and $g(y)=0$. Then the element $x+y$ satisfies (b) (and
(a), (c)). So it suffices to find $x\in G_\delta$ satisfying (a)+(c). If this
fails then for every $x_1, x_2\in G_\delta$ we have
\[[\bigwedge_{h\in R} h(x_1)=h(x_2)]\ \Rightarrow\ g(x_1)=g(x_2).\]
So there is a function $F:{}^R\bbz\longrightarrow\bbz$ such that
$g(x)=F(\dots,h(x),\dots)_{h\in R}$. Take $i<\Gz(M^1_\delta)$ such that
$R\in N^\delta_i$. Then $N^\delta_{i+1}\models \|R\|\leq\|\delta\|$
(remember \ref{principle}({\bf f}, {\bf g})) and necessarily $N^\delta_{i+1}
\models \|R\|<\|\delta\|$ (as $\|R\|<\|\delta\|$ in ${\bf V}$). Applying
\ref{principle}({\bf h}) we get that $Q^{M^1_\delta}_4\in N^\delta_{i+1}$
and therefore $N^\delta_{i+1}\models\|R\|=\kappa_0$, where $\kappa_0=\|R\|$
(in {\bf V}). Let $\kappa_1=2^{\kappa_0}$ (so $\kappa_1<\|\delta\|$). Then,
by \ref{cl1}, we get ${\mathcal P}(\kappa_0)\subseteq N^\delta_{i+1}$ and
$N^\delta_{i+1}\models \|{}^{\textstyle\kappa_0}\bbz\times\bbz\|=\kappa_1$. 
Again by \ref{cl1}, we get ${\mathcal P}(\kappa_1)\subseteq N^\delta_{i+1}$. 
But this implies that ${\mathcal P}({}^{\kappa_0}\bbz\times\bbz)\subseteq
N^\delta_{i+1}$ and ${\mathcal P}({}^R \bbz\times\bbz)\subseteq N^\delta_{i
+1}$. In particular $F\in N^\delta$. Since $\{G_\delta, R\}\in N^\delta$
($G_\delta$ by clause (h) of \ref{principle}, $R$ by the assumption) we
conclude that $g\in N^\delta$ --- a contradiction to condition (h) of
\ref{principle}. The case when $R$ is finite is much easier. We start as
above, but getting $x\in G_\delta$ with (a)+(c) we give purely algebraical
arguments. The claim is proved. 
\end{proof}

\begin{claim}
\label{cl2A}
\begin{enumerate}
\item Suppose that $\cf(\delta)\neq\cf(\mu)<\mu<\delta<\mu^+$ (so we are
in the case $\lambda=\mu^+$), $R\in N^\delta_i$, $R\subseteq\HOM(G_\delta,
\bbz)\cup\bigcup\limits_{q\in\bP}\HOM(G_\delta,\bbz/q\bbz)$, $i<\Gz(
M^1_\delta)$ and $N^\delta_i\models\|R\|=\mu$. Then for each sufficiently
large $\beta\in(\mu,\delta)\setminus S$ there exist $x_j\in G_{\beta+1}$ for
$j<\cf(\mu)$ such that: 
\begin{enumerate}
\item[(a)] if $h\in R$ then the set $\{j<\cf(\mu):h(x_j)\neq 0\}$ is bounded
in $\cf(\mu)$,
\item[(b)] $G_\beta\oplus\lk x_j:j<\cf(\mu)\rk_{G_{\beta+1}}$ is a direct
summand of $G_{\beta+1}$ (and so of $G_\delta$),
\item[(c)] $g(x_j)\neq 0$ for all $j<\cf(\mu)$.
\end{enumerate}
\item In 1), if we change the assumptions to:
\[\cf(\delta)=\theta^+=\mu<\delta<\mu^+,\quad\quad \|R\|=\theta\]
then the assertion holds true after replacing $\cf(\mu)$ by $\cf(\theta)$
(so $x_j$ are being chosen for $j<\cf(\theta)$).
\item Suppose that $\cf(\delta)\neq\cf(\theta)$, $\cf(\delta)<\theta^+=\mu<
\delta<\mu^+$, $R\in N^\delta_i$, $R\subseteq\HOM(G_\delta,\bbz)\cup
\bigcup\limits_{q\in\bP}\HOM(G_\delta,\bbz/q\bbz)$, $i<\Gz(M^1_\delta)$ and 
$N^1_\delta\models\|R\|=\mu$. Then for sufficiently large $\beta\in(\mu,
\delta)\setminus S$ there exists a sequence $\lk x_{j,k}: j<\mu, k<\cf(
\theta)\rk \subseteq G_{\beta+1}$ such that
\begin{enumerate}
\item[(a)] if $h\in R$ then $\{(j,k)\in\mu\times\cf(\theta): h(x_{j,k})\neq
0\}\in I_{\lk\mu,\cf(\theta)\rk}$,
\item[(b)] $G_\beta\oplus\lk x_{j,k}:j<\mu, k<\cf(\theta)\rk_{G_{\beta+1}}$
is a direct summand of $G_{\beta+1}$,
\item[(c)] $g(x_{j,k})\neq 0$ for all $j<\mu$, $k<\cf(\theta)$.
\end{enumerate}
\end{enumerate}
\end{claim}

\begin{proof}[Proof of the claim]
We follow exactly the lines of the proof of \ref{cl2}, but first we have to
choose an increasing sequence $\lk R_j: j<\cf(\mu)\rk\in N^\delta_i$ such
that $\bigcup\limits_{j<\cf(\mu)} R_j=R$, $\|R_j\|<\|R\|$ (and hence
$2^{\|R_j\|}<\mu$ as $\mu$ is a limit cardinal). To find the $R_j$ use
condition (14) of the construction (and $Q_3$, $Q_4$). Then use \ref{cl2} to
find $\beta_0\in (\mu,\delta)$ such that there are $x^*_j\in G_{\beta_0}$
(for $j<\cf(\mu)$) with
\[(\forall h\in R_j)(h(x^*_j)=0)\ \mbox{ and }\ g(x^*_j)\neq 0\]
(remember that $\cf(\delta)\neq\cf(\mu)$). Now if $\beta\in (\beta_0,\delta)
\setminus S$ then we find a sequence $\langle y_j:j<\cf(\mu)\rangle\subseteq
G_{\beta+1}$ which can be extended to a basis of $G_{\beta+1}$ over
$G_\beta$ and such that for all $j<\cf(\mu)$ 
\[(\forall h\in R_j)(h(y_j)=0)\ \mbox{ and }\ g(y_j)=0.\]
Put $x_j=y_j+x^*_j$.
\medskip

\noindent 2)\ \ Similarly (note that if $R\in N^\delta_i$, $\|R\|=\theta$
then $N^\delta_i\models \|R\|=\theta$).
\medskip

\noindent 3)\ \ Similarly: first find $\lk R_{j,k}: j<\mu,k<\cf(\theta)\rk
\in N^\delta_i$ such that $\|R_{j,k}\|<\theta$, the sequence $\lk
\bigcup\limits_{k<\cf(\theta)}R_{j,k}: j<\mu\rk$ is increasing, for each
$j<\mu$ the sequence $\lk R_{j,k}: k<\cf(\theta)\rk$ is increasing and 
$\bigcup\limits_{j<\mu}\bigcup\limits_{k<\cf(\theta)}R_{j,k}= R$. Next
follow as in 1).
\end{proof}

Now we are going to finish the proof of $(\otimes)_4$ (in the case
$\Gz(M^1_\delta)\neq 0$). Note that by $(\otimes)_3$, $f^{p,\zeta^\ell}
\rest G_\delta=f^{p,\zeta^\ell}_\delta$, so we can try to apply condition
(12). But condition (12) says ``help only those who can help themselves''.
More specifically we have to prove that there are $H,f',f^{p,\zeta}_{
\delta+1}$ as required there (in particular $(*)$) and then by (12), $g\rest
G_{\delta+1}$ gives the desired contradiction. But this is done by the
following claim. 

\begin{claim}
\label{cl4}
Suppose that $\delta$, $g\ldots$ are as chosen earlier. Then there exists a
free group $H$ such that $G_\delta\subseteq H$, $H/G_\beta$ is free for
$\beta\in\delta\setminus S$, $\|H\|=\|\delta\|$, the homomorphism $g$ cannot
be extended to a homomorphism $g'\in\HOM(H,\bbz)$ and
\begin{enumerate}
\item[$(\alpha)$] for every $h\in\HOM(G_\delta,\bbz)\cap N^\delta$ and
$h^+\in\HOM(H,\bbz/p\bbz)$ such that $h/p\bbz\subseteq h^+$ there is $h^*
\in\HOM(H,\bbz)$ with $h\subseteq h^*$ and $h^*/p\bbz=h^+$;
\item[$(\beta)$]  if $q\in\bP$, $\zeta<\lambda_q\cap\delta$ then
$f^{q,\zeta}_\delta$ can be extended to an element of $\HOM(H,\bbz/q\bbz)$. 
\end{enumerate}
\end{claim}

\begin{proof}[Proof of the claim]
Let us recall that $\Gz(M^1_\delta)=\cf(\delta)$. For $\varepsilon<\Gz
(M^1_\delta)$ let $R_\vare=\GB^\delta_\varepsilon (M^1_\delta)\cap \big(
\HOM(G_\delta,\bbz)\cup\bigcup\limits_{q\in\bP}\HOM(G_\delta,\bbz/q\bbz)
\big)$ and let $\langle f_\vare: \vare<\Gz(M^1_\delta)\rangle$ 
be a sequence of functions such that for each $\vare<\Gz(M^1_\delta)$
we have:
\[f_\vare:\delta\stackrel{\rm onto}{\longrightarrow} R_\vare\ \ \mbox{ and
}\ \ \langle f_\zeta: \zeta\leq\vare\rangle\in N^\delta_{\vare+1}\]
(see clauses {\bf (e)} and {\bf (g)} of \ref{principle}). Let $\lk
\beta^\delta_i:i<\cf(\delta)\rk$ be an increasing sequence with limit
$\delta$ such that $\beta^\delta_0>\cf(\delta)$ and $\lk\beta^\delta_i:i<
i^*\rk\in N^\delta$ for all $i^*<\cf(\delta)$ (see claim \ref{cl1A}). 
Finally for $i<\cf(\delta)$ let 
\[R^*_i\stackrel{\rm def}{=}\{f_\vare(\zeta): \zeta<\beta^\delta_i\ \&\
\vare\leq i\}.\]
Then $R^*_i$ are increasing with $i$ and
\[\bigcup_{i<\cf(\delta)} R^*_i = \bigcup_{\vare<\Gz(M^1_\delta)}
R_\vare\] 
and for each $i^*<\cf(\delta)$ the sequence $\lk R^*_i: i<i^*\rk$ belongs to
$N^\delta_{\zeta_{i^*}}$ (for some $\zeta_{i^*}<\Gz(M^1_\delta)$). Moreover
$N^\delta_{\zeta_{i^*}}\models\|R^*_i\|=\|\beta^\delta_i\|$ for each $i<i^*$.

Let $\lk \alpha^\delta_i: i<\cf(\delta)\rk$ be an increasing continuous
sequence cofinal in $\delta$ and disjoint from $S$ (possible by the choice
of $S$).  
\medskip

\begin{description}
\item[\bf Case A] $\delta$ is a strongly limit singular cardinal.
\end{description}
In this case we have
\[(\forall i<\cf(\delta))(2^{\|R^*_i\|}<\delta).\]
Thus we may apply claim \ref{cl2} and choose by induction on
$i<\cf(\delta)$ an increasing sequence $\lk j_i: i<\cf(\delta)\rk\subseteq
\cf(\delta)$ and $x^\delta_i\in G_{\alpha^\delta_{j_i}+1}$ such that:
\begin{enumerate}
\item[(a)] $h\in R^*_i\Rightarrow h(x^\delta_i)=0$;
\item[(b)] $G_{\alpha^\delta_{j_i}}\oplus({\mathbb Z}x^\delta_i)$ is a direct
summand of $G_\delta$;
\item[(c)] $g(x^\delta_i)\neq 0$.
\end{enumerate}
Since $\lk \alpha^\delta_i: i<\cf(\delta)\rk \subseteq\delta\setminus S$ is
increasing continuous (and cofinal in $\delta$) we get that the subgroup
$H_\delta=\lk x^\delta_i: i<\cf(\delta)\rk_{G_\delta}$ is a direct summand
of $G_\delta$, say $G_\delta=H_\delta\oplus H^*_\delta$ (and $\{x^\delta_i:
i<\cf(\delta)\}$ is a free basis of $H_\delta$).  Let $I=\{A \subseteq
\cf(\delta): A$ is bounded$\}$ and apply \ref{1.2B} for $\cf(\delta)$, $I$
and $H_\delta$ and get the respective free group $H'\supseteq H_\delta$
($H'\cap G_\delta=H_\delta$). We claim that the group $H=H'\oplus
H^*_\delta$ is as required. For this first note that if $\beta\in\delta
\setminus S$, $\alpha^\delta_{j_{i_0}}>\beta$, $A=[0,j_{i_0})$ then $H'/\lk
x^\delta_i:i\in A\rk$ is free and hence $H/G_{\alpha^\delta_{j_{i_0}}}$ is
free. But $G_{\alpha^\delta_{j_{i_0}}}/G_\beta$ is free so we conclude that
$H/G_\beta$ is free.

As $g$ cannot be extended to a member of $\HOM(H',\bbz)$ it has no extension
in $\HOM(H,\bbz)$. Suppose now that $h\in\HOM(G_\delta,\bbz)\cap N^\delta$,
so $h\in R^*_{i_0}$ for some $i_0<\cf(\delta)$. Let $h^+\in\HOM(H,\bbz/
p\bbz)$ extend $h/p\bbz$. Since for all $i\geq i_0$ we have $h(x^\delta_i)
=0$, we may apply clause $(\gamma)$ of \ref{1.2B} to get a suitable lifting
$h^*\in\HOM(H,\bbz)$ of $h^+$. Similarly, we use \ref{1.2B}($\delta$) to
show that $f^{q,\xi}_\delta$ can be extended onto $H$ (remember
$f^{q,\xi}_\delta\in N^\delta$). 
\medskip

\begin{description}
\item[\bf Case B] $\delta\in (\mu,\mu^+)$ for some cardinal $\mu$
such that $\cf(\mu)=\cf(\delta)<\mu$.
\end{description}
By condition (14) of the construction and the use of $Q_3$,
$Q_4$ we have that, letting $\alpha=\cf(\delta)$, for each $i<\cf(\delta)$
\[N^\delta_{\zeta_{i+1}}\models\mbox{`` }\|\beta^\delta_i\|=\|\mu\|\ \ \&\ \
\cf(\mu)=\alpha\mbox{ ''.}\] 
This allows us to build $R^{**}_i$ such that 
\[i<j<\cf(\delta)\ \ \Rightarrow\ \ R^{**}_i\subseteq R^{**}_j\in
N^\delta,\]
\[\bigcup_{i<\cf(\delta)}R^{**}_i=\bigcup_{i<\cf(\delta)}R^*_i\ \ \mbox{ and
}\ \ \|R^{**}_i\|<\mu.\]
Now we can continue as in the previous case.
\medskip

\begin{description}
\item[\bf Case C]  $\delta\in (\mu,\mu^+)$ for some cardinal number $\mu$
such that $\cf(\delta)\neq\cf(\mu)<\mu$
\end{description}
By claim \ref{cl2A}(1) we can choose $x^\delta_{i, j}$ (for $i<\cf(\delta)$,
$j<\cf(\mu)$) such that:
\begin{enumerate}
\item[(a)] for each $h\in R^*_i$, for every $j<\cf(\mu)$ large enough
$h(x^\delta_{i, j} )=0$;
\item[(b)] $\{x^\delta_{i, j}:i<\cf(\delta), j< \cf(\mu)\}$ is a free basis
of a direct summand of $G_\delta$, moreover for some increasing sequence
$\langle j_i: i<\cf(\delta)\rangle\subseteq\cf(\delta)$, for each $i^*<\cf(
\delta)$, the family $\{x^\delta_{i, j}: i^*\le i<\cf(\delta), j<\cf(\mu)\}$
is a free basis of a subgroup $H\subseteq G_\delta$ such that
$G_{\alpha^\delta_{j_{i^*}}}\oplus H$ is a direct summand of $G_\delta$;
\item[(c)] $g(x^\delta_{i, j})\neq 0$.
\end{enumerate}
Let 
\[I=I_{\lk \cf(\delta),\cf(\mu)\rk}=\{A\subseteq\cf(\delta)\times\cf(\mu):
(\forall^* i<\cf(\delta))(\forall^* j<\cf(\mu))((i,j)\notin A)\}.\]
Again apply \ref{1.2B} (with $\lambda$ there standing for $\cf(\delta)+
\cf(\mu)$).
\medskip

\begin{description}
\item[\bf Case D]  $\delta\in (\mu,\mu^+)$ for some cardinal number $\mu$
such that $\cf(\delta)=\mu=\cf(\mu)$ is inaccessible.
\end{description}
Similar to Case B.
\medskip

\begin{description}
\item[\bf Case E]  $\delta\in (\mu,\mu^+)$ for some cardinal number $\mu$
such that $\mu=\theta^+=\cf(\delta)$.
\end{description}
First find an increasing sequence $\lk R^{**}_i: i<\theta^+\rk$ such that
$R^{**}_i\in N^\delta$, $\bigcup\limits_{i<\theta^+}R^{**}_i =
\bigcup\limits_{i<\cf(\delta)} R^*_i$ and $\|R^{**}_i\|\leq\theta$. Then
apply claim \ref{cl2A}(2) to choose a sequence $\lk x^\delta_{i,j}:
i<\cf(\delta), j<\cf(\theta)\rk$ similarly as in case C.  
\medskip

\begin{description}
\item[\bf Case F]  $\delta\in (\mu,\mu^+)$ for some cardinal $\mu$
such that $\cf(\delta)<\mu=\theta^+$, $\cf(\delta)\neq\cf(\theta)$.
\end{description}
Using claim \ref{cl2A}(3) we choose an increasing sequence $\langle j_i:i
<\cf(\delta)\rangle\subseteq\cf(\delta)$ and a sequence $\lk x^\delta_{i,j,
k}:i<\cf(\delta), j<\mu,k<\cf(\theta)\rk$ such that
\begin{enumerate}
\item[(a)] for each $h\in R_i^*$, for every $j<\mu$ large enough for every
$\kappa<\cf(\theta)$ large enough, $h(x^\delta_{i, j, k})=0$;
\item[(b)] $\{x^\delta_{i,j,k}:i<\cf(\delta), j<\mu,\kappa<\cf(\theta)\}$
is a free basis of a direct summand of $G_\delta$; moreover for each
$i^*<\cf(\delta)$ the set 
\[\{ x^\delta_{i,j,k}: i^*\le i<\cf(\delta), j<\mu,k<\cf(\theta)\}\]
is a free basis of a subgroup $H\subseteq G_\delta$ such that
$G_{\alpha^\delta_{j_{i^*}}}\oplus  H$ is a direct summand of $G_\delta$;
\item[(c)] $g(x^\delta_{i,j,k})\neq 0$.
\end{enumerate}
Let $I=I_{\lk \cf(\delta),\mu,\cf(\theta)\rk}$ and apply \ref{1.2B}.
\medskip

\begin{description}
\item[\bf Case G]  $\delta\in (\mu,\mu^+)$ for some cardinal $\mu$
such that $\mu=\theta^+$, $\cf(\delta)=\cf(\theta)$.
\end{description}
This is similar to the case F though we have to modify the application of
\ref{1.2B}. First we choose increasing sequences $\lk R^{**}_{i,j}:j<\mu\rk
\in N^\delta$ (for $i<\cf(\delta)$) such that 
\[\|R^{**}_{i,j}\|<\mu,\quad\quad\bigcup_{j<\mu}
R^{**}_{i,j}=R^*_i\quad\mbox{ and}\]
\[\lk\lk R^{**}_{i,j}:j<\mu\rk: i<i^*\rk\in N^\delta\quad\quad\mbox{ for
each } i^*<\cf(\delta).\]
Then we choose increasing sequences $\lk R^{***}_{i,j,k}: k<\cf(\theta)\rk$
(for $i<\cf(\delta)=\cf(\theta)$, $j<\mu$) such that 
\[\|R^{***}_{i,j,k}\|<\theta,\quad\quad\bigcup_{k<\cf(\theta)}
R^{***}_{i,j,k} = R^{**}_{i,j}\quad\mbox{ and}\]
\[\lk\lk R^{***}_{i,j,k}: j<\mu, k<\cf(\theta)\rk: i<i^*\rk\in N^\delta.\]
Now, for $\ell<\cf(\theta)=\cf(\delta)$, $j<\mu$ put
$R^{+}_{j,\ell}=\bigcup\limits_{i,k<\ell} R^{***}_{i,j,k}$. Note that
$R^+_{j,\ell}\in N^\delta$ and $\|R^+_{j,\ell}\|<\theta$. Moreover if
$h\in\bigcup\limits_{i<\cf(\delta)}R^*_i$ then $(\forall^* j<\mu)(\forall^*
\ell<\cf(\delta))(h\in R^+_{j,\ell})$. Next, as in the proof of \ref{cl2A}
we choose $x^*_{j,\ell}$, $y_{j,\ell}$ such that 
\smallskip

$(\forall h\in R^+_{j,\ell})(h(x^*_{j,\ell})=h(y_{j,\ell})=0)$,\quad $g(
x^*_{j,\ell})\neq 0$, $g(y_{j,\ell})=0$,

if $\rho_{j,\ell}\stackrel{\rm def}{=}\min\{\alpha^\delta_{\ell_0}:
x^*_{j,\ell}\in G_{\alpha^\delta_{\ell_0}}\}$ then $\ell<\ell_0$ and

$\{y_{j,\ell}: \rho_{j,\ell}=\beta\}\subseteq G_{\beta+1}$ can be extended
to a basis of $G_{\beta+1}$ over $G_\beta$.
\smallskip

\noindent Then we put $x_{j,\ell}=x^*_{j,\ell}+y_{j,\ell}$ (for $j<\mu$,
$\ell<\cf(\delta)$) and we apply \ref{1.2B} as earlier.
\medskip

\begin{description}
\item[\bf Case H]  $\delta\in (\mu,\mu^+)$ for some inaccessible cardinal
$\mu$ such that $\cf(\delta)<\cf(\mu)=\mu$.
\end{description}
This is similar to case C.
\medskip

This completes the proof of claim \ref{cl4}.\end{proof}

The case $\Gz(M^1_\delta)=0$ is much easier and can be done similarly. We do
not have $N^\delta$ and we have to take care of extending homomorphisms
$f^{q,\zeta}_\delta$ only. We basically follow the lines of the previous
case, but proving the suitable variants of \ref{cl2}, \ref{cl2A} instead of
the fact that $g\notin N^\delta$ we use clause (6) of the inductive
construction. 

This completes the proof of $(\otimes)_4$. 
\medskip

To finish the proof of the theorem we have to show 
\begin{enumerate}
\item[$(\otimes)_5$] if $p\in\bP$, $f\in\HOM(G,\bbz/p\bbz)$ then there are
$n<\omega$, $\zeta^0,\ldots,\zeta^{n-1}<\lambda_p$,
$a_0,\ldots,a_{n-1}\in\bbz/p\bbz$ such that $f-\sum\limits_{\ell<n}a_\ell
f^{p,\zeta^\ell}\in\HOM^-(G,\bbz/p\bbz)$ (i.e. the difference can be lifted
to a homomorphism to $\bbz$). 
\end{enumerate}
For this we inductively define a sequence $\lk f_\xi:\xi<\xi(*)\rk\subseteq
\HOM(G,\bbz/p\bbz)$ by:
\begin{quotation}
\noindent $f_\xi$ is the $<^*$-first member of the vector space
$\HOM(G,\bbz/p\bbz)$ over the field $\bbz/p\bbz$ which does not depend on
\[\{f^{p,\zeta}:\zeta<\lambda_p\}\cup\{f_\zeta: \zeta<\xi\}.\]
\end{quotation}
(So $\xi(*)$ is the maximal length of a sequence with the property stated
above.) It is enough to show that all the homomorphisms $f_\xi$ (for
$\xi<\xi(*)$) can be lifted to homomorphisms $f^*_\xi\in\HOM(G,\bbz)$. But
by $(\otimes)_4$ we know that $M^1_\lambda\in\bk_\lambda\setminus
\bk^-_\lambda$, so we may apply {\bf (j)}. Thus we have a club $C\subseteq
E$ such that for each $\delta\in C$:
\[\{C\cap\delta, \lk f_\zeta\rest\delta: \zeta<\delta\rk\}\in N^\delta.\]
Applying the ``moreover'' part of {\bf (j)} of \ref{principle} we may make
use of conditions (8), (9) of the construction. (Remember that in (7) the
sequence $\lk f^{p,\xi}_\alpha: \lambda_p\cap\alpha<\xi<\zeta(p,\alpha)\rk$
has the same definition as our sequence $\lk f_\xi:\xi<\xi(*)\rk$.) 
\end{proof}

\begin{remark} The main theorem \ref{1.2} can be proved for all regular
cardinals which are not weakly compact. This requires some changes in the
construction (and \ref{1.2B}). 
\end{remark}

\end{document}